%% file: earthquakes.tex
\documentclass[a4paper]{amsart}

\usepackage{amssymb,pstricks,amscd,epsfig}
\usepackage[totalwidth=17cm,totalheight=24cm]{geometry}
\usepackage{graphicx}
\input{diagrams.tex}
\diagramstyle[labelstyle=\scriptstyle]

\begin{document}
\input{macros}

\title{Fixed points of compositions of earthquakes}
\author{Francesco Bonsante}
\address{Universit\`a degli Studi di Pavia\\
Via Ferrata, 1\\
27100 Pavia, Italy}
\email{bonsante@sns.it}
\author{Jean-Marc Schlenker}
\thanks{J.-M. S. was partially supported by the A.N.R. through projects
RepSurfaces, ANR-06-BLAN-0311, and GeomEinstein, ANR-09-BLAN-0116-01.}
\thanks{F.B. is partially supported by the A.N.R. through project Geodycos }
\address{Institut de Math\'ematiques de Toulouse, UMR CNRS 5219 \\
Universit\'e Toulouse III \\
31062 Toulouse cedex 9, France}
\email{schlenker@math.univ-toulouse.fr}
\date{April 2010 (v1)}

\begin{abstract}
Let $S$ be a closed surface of genus at least $2$, and let $\lambda$ and $\mu$
be two laminations that fill $S$. Let $E_r^\lambda$ and $E_r^\mu$ be the right
earthquakes on $\lambda$ and $\mu$ respectively. We show that the composition 
$E_r^\lambda\circ E_r^\mu$ has a fixed point in the Teichm\"uller space of $S$.
Another way to state this result is that it is possible to prescribe any two
measured laminations that fill a surfaces as the upper and lower measured 
bending laminations of the convex core of a globally hyperbolic AdS manifold. 
The proof uses some estimates from the geometry of those AdS manifolds.
\end{abstract}

\maketitle


\input jms1 

\input fb1 

\input jms2 

\input fb3 

\appendix

\input recurrence

\input fb2 

\bibliographystyle{amsplain}
\bibliography{/home/schlenker/papiers/outils/biblio}

\end{document}

%% file: macros.tex
\newtheorem{cor}{Corollary}[section]
\newtheorem{theorem}[cor]{Theorem}
\newtheorem{prop}[cor]{Proposition}
\newtheorem{lemma}[cor]{Lemma}
\newtheorem{sublemma}[cor]{Sublemma}
\theoremstyle{definition}
\newtheorem{defi}[cor]{Definition}
\theoremstyle{remark}
\newtheorem{remark}[cor]{Remark}
\newtheorem{example}[cor]{Example}

\newcommand{\cD}{{\mathcal D}}
\newcommand{\cF}{{\mathcal F}}
\newcommand{\cM}{{\mathcal M}}
\newcommand{\cN}{{\mathcal N}}
\newcommand{\cT}{{\mathcal T}}
\newcommand{\cML}{{\mathcal M\mathcal L}}
\newcommand{\cFML}{{\mathcal F\mathcal M\mathcal L}}
\newcommand{\cGH}{{\mathcal G\mathcal H}}
\newcommand{\cQF}{{\mathcal Q\mathcal F}}
\newcommand{\C}{{\mathbb C}}
\newcommand{\N}{{\mathbb N}}
\newcommand{\R}{{\mathbb R}}
\newcommand{\Z}{{\mathbb Z}}
\newcommand{\Kt}{\tilde{K}}
\newcommand{\Mt}{\tilde{M}}
\newcommand{\dr}{{\partial}}
\newcommand{\kappab}{\overline{\kappa}}
\newcommand{\pib}{\overline{\pi}}
\newcommand{\taub}{\overline{\tau}}
\newcommand{\ub}{\overline{u}}
\newcommand{\Sigmab}{\overline{\Sigma}}
\newcommand{\gd}{\dot{g}}
\newcommand{\diff}{\mbox{Diff}}
\newcommand{\dev}{\mbox{dev}}
\newcommand{\devb}{\overline{\mbox{dev}}}
\newcommand{\devt}{\tilde{\mbox{dev}}}
\newcommand{\vol}{\mbox{Vol}}
\newcommand{\hess}{\mbox{Hess}}
\newcommand{\cb}{\overline{c}}
\newcommand{\db}{\overline{\partial}}
\newcommand{\Sigmat}{\tilde{\Sigma}}

\newcommand{\cunc}{{\mathcal C}^\infty_c}
\newcommand{\cun}{{\mathcal C}^\infty}
\newcommand{\dd}{d_D}
\newcommand{\dmin}{d_{\mathrm{min}}}
\newcommand{\dmax}{d_{\mathrm{max}}}
\newcommand{\Dom}{\mathrm{Dom}}
\newcommand{\dn}{d_\nabla}
\newcommand{\ded}{\delta_D}
\newcommand{\delmin}{\delta_{\mathrm{min}}}
\newcommand{\delmax}{\delta_{\mathrm{max}}}
\newcommand{\hmin}{H_{\mathrm{min}}}
\newcommand{\maxi}{\mathrm{max}}
\newcommand{\oL}{\overline{L}}
\newcommand{\oP}{{\overline{P}}}
\newcommand{\xb}{{\overline{x}}}
\newcommand{\yb}{{\overline{y}}}
\newcommand{\Ran}{\mathrm{Ran}}
\newcommand{\tgamma}{\tilde{\gamma}}
\newcommand{\cotan}{\mbox{cotan}}
\newcommand{\area}{\mbox{Area}}
\newcommand{\lambdat}{\tilde\lambda}
\newcommand{\xt}{\tilde x}
\newcommand{\Ct}{\tilde C}
\newcommand{\St}{\tilde S}

\newcommand{\sh}{\mathrm{sinh}\,}
\newcommand{\ch}{\mathrm{cosh}\,}

\newcommand{\II}{I\hspace{-0.1cm}I}
\newcommand{\III}{I\hspace{-0.1cm}I\hspace{-0.1cm}I}
\newcommand{\note}[1]{{\small {\color[rgb]{1,0,0} #1}}}

%% file: jms1.tex
\section{Introduction, main results}

\subsection{Earthquakes.}

In this paper we consider a closed surface $S$ of genus at least $2$.
We denote by $\cT_S$, or sometimes simply by $\cT$, the Teichm\"uller 
space of $S$, which we consider to be the space of hyperbolic metrics
on $S$ considered up to isotopy. We denote by $\cML_S$, or simply by
$\cML$, the space of measured laminations on $S$.

Given a measured lamination $\lambda\in \cML_S$, we denote by
$E_l^\lambda$ the {\it left earthquake} along $\lambda$ on $S$.
$E_l^\lambda$ is a continuous map from $\cT_S$ to $\cT_S$, see
e.g. \cite{thurston-notes, kerckhoff}. Given two measured laminations
$\lambda, \mu$ in $S$, we say that $\lambda$ and $\mu$ {\it fill} $S$
if any closed curve $c$ in $S$ which is not homotopically trivial has
non-zero intersection with either $\lambda$ or $\mu$.

The main result of this paper concerns fixed points of compositions
of earthquakes on $S$.

\begin{theorem} \label{tm:1}
Let $\lambda, \mu\in \cML_S$ be two measured laminations which fill 
$S$. Then $E_l^\lambda\circ E_l^\mu:\cT_S\rightarrow \cT_S$ has a fixed point in 
$\cT_S$. 
\end{theorem}

There are some reasons to believe that this fixed point is unique.
It is explained below why this statement is equivalent to a conjecture
made by G. Mess \cite{mess,mess-notes} concerning globally hyperbolic
anti-de Sitter $3$-manifolds.

Theorem \ref{tm:1} shows a contrast between earthquakes and Dehn twists.
Let $\lambda$ and $\mu$ be two closed curves that fill $S$, and let 
$D^\lambda$ and $D^\mu$ be the Dehn twists along $\lambda$ and $\mu$, 
respectively. Thurston proved (see \cite[Section 6]{thurston:surfaces})
that for $n$ and $m$ large enough, the composition $(D^\lambda)^n\circ
(D^\mu)^m$ is a pseudo-Anosov diffeomorphism of $S$, so that it acts on 
$\cT$ without fixed point. This does not contradict Theorem \ref{tm:1}
since Dehn twists do not act on $\cT$ as earthquakes. (The shearing 
distance of a Dehn twist depends on the hyperbolic metric.)

\subsection{Quasifuchsian hyperbolic 3-manifolds.}

Let $M:=S\times \R$, let $\cQF_S$ be the space of quasifuchsian hyperbolic metrics
on $M$, that is, complete hyperbolic metrics containing a non-empty compact subset
which is convex, considered up to isotopy. Note that, here and elsewhere, we 
say that a subset $K\subset M$ is {\it convex} if any geodesic segment 
with endpoints in $K$ is contained in $K$. If
$K$ is convex then its lift to the universal cover of $M$ is also convex and
therefore connected, and, if $K$ is non-empty, it is a deformation retract of $M$.

Quasifuchsian hyperbolic metrics have a boundary at infinity
which is topologically the 
union of two copies of $S$, endowed with a 
complex structure. According to a classical theorem of Bers \cite{bers},
this complex structure defines a parameterization of $\cQF_S$ by the 
product of two copies of $\cT_S$, one for each boundary component.

Given $g\in \cQF_S$, $(M,g)$ contains a
smallest non-empty closed convex subset $C(M)$, called its convex core.
Except when $C(M)$ is a totally geodesic surface, 
its boundary is the disjoint union of two closed surfaces
each homeomorphic to $S$. We call them the ``upper'' and ``lower'' boundary
components of $C(M)$. Each has an induced metric which is hyperbolic, and
is ``bent'' along a measured geodesic lamination, see \cite{thurston-notes,epstein-marden}.

Let $m_+, m_-$ be the induced hyperbolic metrics and $\lambda_+,\lambda_-$
be the measured bending laminations on the upper and lower boundary
components of $C(M)$. Then it is well known that $\lambda_-$ and 
$\lambda_+$ fill $S$. Moreover 
$\lambda_+$ and $\lambda_-$ have no closed curve with weight larger than
$\pi$. 

Thurston conjectured that any two measured laminations on $S$ satisfying
these two conditions can be uniquely realized as $(\lambda_-,\lambda_+)$.
The existence part of this conjecture was proved (in a more general form) 
by Bonahon and Otal.

\begin{theorem}[Bonahon, Otal \cite{bonahon-otal}] \label{tm:bo}
Let $\lambda_+,\lambda_-$ be two measured laminations which fill $S$, 
which have no closed curve with weight at least equal to $\pi$. There exists
a quasifuchsian metric $g\in \cQF_S$ such that the measured bending
lamination on the upper (resp. lower) boundary component of $C(M)$
is $\lambda_+$ (resp. $\lambda_-$).
\end{theorem}

Thurston also conjectured that any two hyperbolic metrics can be obtained
uniquely as $(m_+,m_-)$, and the existence part of this statement is also
known (it follows from \cite{L4,epstein-marden}). This type of statement
however does not appear explicitly here.

\subsection{Globally hyperbolic AdS manifolds.}

The geometric theory of 3-dimensional globally hyperbolic anti-de Sitter
(AdS) manifolds is quite similar to the theory of quasifuchsian hyperbolic
3-manifolds, a remarkable fact discovered by Mess \cite{mess,mess-notes}.

Recall that an AdS 3-manifold $M$ is {\it globally hyperbolic} if it
contains a space-like surface $S$ which intersects any inextendible
time-like curve exactly once. It is {\it globally hyperbolic maximal
  compact} (GHMC) if this space-like surface is closed, and moreover
$M$ is maximal (for the inclusion) under the previous condition.

\begin{remark}
Global hyperbolicity is a notion that makes sense for every Lorentzian
manifold. It has strong consequences. If $S$ is a Cauchy surface in
$M$, then topologically $M=S\times\mathbb R$.

Moreover it can be shown that there is a foliation of $M$ into
spacelike slices parallel to $S$. More precisely, for  a suitable
product structure on $M$, the metric is of the
form
\begin{equation}\label{eq1}
    -dt^2+g_t
\end{equation}
where $t$ is the real parameter and $g_t$ is a path of Riemann metrics
on $S$.

Conversely, assuming $S$ to be compact, any metric on $S\times\mathbb R$ 
that is of the form (\ref{eq1}) is globally hyperbolic.

We refer to \cite{beem, geroch} for a  complete treatment of this topic.
\end{remark}

We call $\cGH_S$ the space of globally hyperbolic
maximal compact AdS metrics on $S\times \R$, with the Cauchy surface homeomorphic
to $S$, considered up to isotopy.

Mess \cite{mess,mess-notes} discovered that to any AdS metric $g\in \cGH_S$ are
associated two points in $\cT_S$, its left and right hyperbolic metric, which can be
defined through the decomposition of the identity component of the isometry
group of $AdS_3$ as the product of two copies of $PSL_2(\R)$. 
Moreover, these left and right metrics
define a parameterization of $\cGH_S$ by $\cT_S\times \cT_S$. This can
be construed as an analog of the Bers theorem mentioned above. 
Given $g\in \cGH_S$, $(M,g)$ also contains a smallest closed non-empty convex 
subset which is compact, which we also call $C(M)$ here. 

We say that $(M,g)$ is {\it Fuchsian} if $C(M)$ is a totally geodesic surface.
As in the quasifuchsian 
setting, if $M$ is not Fuchsian, then the boundary of $C(M)$
is the union of two surfaces homeomorphic to $S$, and
each has a hyperbolic induced metric and is bent along a measured geodesic
lamination. We call $m_+, m_-$ the two hyperbolic metrics, and $\lambda_+,
\lambda_-$ the two measured laminations. Extending Thurston's conjectures,
Mess \cite{mess} asked whether any two hyperbolic metrics on $S$ can be 
uniquely obtained, and whether any two measured laminations that fill $S$
can be uniquely obtained. The second result presented here is the proof
of the existence part of the statement concerning the measured bending
laminations.

\begin{theorem}\label{tm:2}
Let $\lambda_+,\lambda_-\in \cML_S$ be two measured laminations that
fill $S$. There exists a globally hyperbolic maximal AdS manifold
$M$ such that $\lambda_+$ and $\lambda_-$ are the measured pleating
laminations on the upper and lower boundary components of the 
convex core. 
\end{theorem}

\subsection{From earthquakes to AdS manifolds.}

Theorem \ref{tm:2} is equivalent to Theorem \ref{tm:1} thanks to some
key properties of GHMC AdS manifolds, which we recall briefly here and
in more details in Section 2. 

\begin{theorem}[Mess \cite{mess}] \label{tm:mess}
\begin{itemize}
\item Given a GHMC AdS manifold, its left and right metrics 
$\rho_l,\rho_r$ are related to the induced metrics and measured 
bending laminations on the boundary of the convex core as follows:
$$ \rho_l = E_l^{\lambda_+}(m_+)~, ~~ \rho_r = E_r^{\lambda_+}(m_+)~, $$
$$ \rho_l = E_r^{\lambda_-}(m_-)~, ~~ \rho_r = E_l^{\lambda_-}(m_-)~, $$
so that 
\begin{equation}\label{eq:mess}
\rho_l = E_l^{2\lambda_+}(\rho_r) = E_r^{2\lambda_-}(\rho_r)~. 
\end{equation}
\item Given $\rho_l,\rho_r\in \cT_S$ and $\lambda_+,\lambda_-\in \cML$
such that Equation (\ref{eq:mess}) holds, there is a unique GHMC AdS
manifold with left and right metrics $\rho_l$ and $\rho_r$,
and the measured bending laminations on the boundary components
of its convex core are $\lambda_+,\lambda_-$.
\end{itemize}
\end{theorem}

The proof of Theorem \ref{tm:2} clearly follows from this and from 
Theorem \ref{tm:1} (and conversely): given two measured laminations
$\lambda_+,\lambda_-$ that fill $S$, the map $E_r^{2\lambda_+}\circ
E_r^{2\lambda_-}$ has a fixed point, which we call $\rho_r$. Setting
$\rho_l:=E_r^{2\lambda_-}(\rho_r)=E_l^{2\lambda_+}(\rho_r)$, we see with
Theorem \ref{tm:mess} that $\rho_l,\rho_r$ are the left and right
metrics of a GHMC AdS manifold for which $\lambda_+$ and
$\lambda_-$ are the upper and lower measured bending laminations of
the boundary of the convex core. The same argument can be used to 
prove Theorem \ref{tm:1} from Theorem \ref{tm:2}. 

\subsection{Outline of the proof.}

The proof of Theorem \ref{tm:2} has two main parts. 

The first is a 
description of elements of $\cGH_S$ near the ``Fuchsian locus'',
that, is the subset of AdS metrics on $S\times \R$ containing a
totally geodesic Cauchy surface. In the quasifuchsian case, 
Bonahon \cite{bonahon-almost} proved that any two measured laminations
which fill $S$ and are ``small enough'' in a suitable sense can be
uniquely obtained as the measured bending laminations on the 
boundary of the convex core of a quasifuchsian manifold which is
``almost-Fuchsian''. Series \cite{series-limits} then proved that
those almost-Fuchsian metrics are the only ones realizing 
$\lambda_+,\lambda_-$ as the bending lamination on the boundary of
the convex core. 

We present here an analog of those arguments for the GHMC AdS setting.

\begin{theorem} \label{tm:4}
Let $\lambda,\mu\in \cML_S$ be two measured laminations that fill $S$.
There exists $\epsilon>0$ such that, for all $t\in (0,\epsilon)$, 
there exists a unique GHMC AdS manifold such that the measured bending
laminations on the upper and lower boundary components of the convex
core are $t\lambda$ and $t\mu$.  
\end{theorem}

The second tool of the proof of Theorem \ref{tm:2} is a compactness
statement. First a definition.

\begin{defi}
Let $\cFML_S\subset\cML_S\times \cML_S$ be the space of pairs of 
measured laminations that fill $S$. Let $\Phi:\cT_S\times
\cT_S\rightarrow \cFML_S$ be the map which associates to $(\rho_l,\rho_r)$
the measured bending laminations $(\lambda_+,\lambda_-)$ on the 
boundary of the convex core of the unique GHMC AdS manifold with
left and right metrics $\rho_l$ and $\rho_r$. 
\end{defi}

The compactness statement that is needed is the following statement,
which is equivalent to Proposition \ref{pr:proper}.

\begin{prop} \label{pr:proper-1}
$\Phi$ is proper.
\end{prop}

Since $\Phi$ is continuous, it is possible to define its degree,
which by Theorem \ref{tm:4} is equal to $1$. Therefore $\Phi$ is
surjective, which is another way to state Theorem \ref{tm:2}.

The proof of Proposition \ref{pr:proper-1} is based on Proposition
\ref{pr:main}, a simple estimate on the intersection between
measured laminations under some geometric assumptions.

\begin{remark}
The degree argument used here is possibly not limited to the AdS
setting. It might also be used for hyperbolic quasifuchsian
3-manifolds to prove Theorem \ref{tm:bo}. The proof should then use
the compactness result of \cite{bonahon-otal} to prove that the map
sending a quasifuchsian metric to the measured laminations on the
boundary of the convex core (as seen as a map to the space of
``admissible'' pairs of laminations) is proper, so that its degree can
be considered. Then the analysis made by Bonahon \cite{bonahon-almost}
concerning the behavior of this map near the Fuchsian locus, and the
results of Series \cite{series-limits} showing that small laminations
can be obtained only there, should indicate that the degree is one, so
that the map is surjective. An interesting facet of this possible
argument would be that it does not use the rigidity of hyperbolic
cone-manifolds proved by Hodgson and Kerckhoff \cite{HK}.
\end{remark}

\subsection{Cone singularities.}

The arguments used here can be extended, with limited efforts, to 
hyperbolic surfaces with cone singularities. Recall that Thurston's
Earthquake Theorem can be extended to this case \cite{cone}. 
We need some notations before stating the result. 

We call $\cT_{S,\theta}$ the Teichm\"uller
space of conformal structures on $S$ with $N$ marked points
$x_1,\cdots, x_N$ of cone angles given by 
$\theta=(\theta_1,\cdots,\theta_N)\in (0,\pi)^N$. Any conformal 
class in $\cT_{S,\theta}$ contains a unique hyperbolic metric with cone
singularities of angle $\theta_i$ at each $x_i$ (see \cite{troyanov}). 
Let $\cML_{S,N}$ be the space of measured laminations on
the complement of $\{ x_1,\cdots, x_N\}$ in $S$. Given a hyperbolic
metric $g$ on $S$ with cone singularities of angle $\theta_i$ at each
$x_i$, any $\lambda\in \cML_{S,N}$ can be realized uniquely as a
measured geodesic lamination for $g$ (as long as the $\theta_i$
are less than $\pi$), see e.g. \cite{dryden-parlier}. 
The notion of earthquake as defined by
Thurston extends to this setting with cone singularities, see
\cite{cone}. In this context, we (still) say that two measured
laminations $\lambda,\mu\in \cML_{S,N}$ {\it fill} $S$ if, given
any closed curve $c$ in $S\setminus \{ x_1,\cdots, x_N\}$ not
homotopic to zero or to one of the $x_i$,
the intersection of $c$ with either $\lambda$ or $\mu$
is non-zero.

\begin{theorem} \label{tm:5}
Let $\lambda,\mu\in \cML_{S,N}$ be two measured laminations which 
fill $S$. Then $E_r^{\lambda}\circ E_r^{\mu}$ has at least one
fixed point on $\cT_{S,\theta}$. 
\end{theorem}

The proof follows the same line as the proof of Theorem \ref{tm:1}.
Theorem \ref{tm:mess} still holds in the context of hyperbolic
surfaces with cone singularities (of fixed angle less than $\pi$), in
this case the GHMC AdS manifolds that are involved have ``particles'',
that is, cone singularities along time-like curves, see
\cite{cone}. Those manifolds have a convex core, its boundary is a
pleated surface (outside the cone singularities) with an induced
metric which is hyperbolic with cone singularities, and the pleating
defines a measured lamination.  Theorem \ref{tm:5} can therefore be
stated in an equivalent way involving measured laminations on the
boundary of the convex core.

\begin{theorem} \label{tm:6}
Let $\lambda_+,\lambda_-\in \cML_{S,N}$ be two measured laminations which fill
$S$ (considered as a surface with $N$ punctures). There exists a GHMC AdS 
manifold with particles, with angles $\theta_1,\cdots, \theta_N$, such that 
the measured bending laminations on the upper and lower boundary components
of the convex core are $\lambda_+$ and $\lambda_-$.
\end{theorem}

The proof of Theorem \ref{tm:6} basically follows the same path as the proof
of Theorem \ref{tm:2}, with ``particles'' in the GHMC AdS manifolds which is
considered.

\subsection{Flat space-times.}

The Minkowski globally hyperbolic spacetimes 
have interesting properties which are reminiscent of,
but different from, those of globally hyperbolic AdS manifolds. We recall
those properties briefly here, more details are given in Appendix B.

Consider again a closed surface $S$ of genus at least $2$. 
Maximal globally hyperbolic spacetimes containing a space-like surface
homeomorphic to $S$ are quotient of a convex domain in the
3-dimensional Minkowski space $\R^{2,1}$ by $\rho(\pi_1(S))$, where
$\rho$ is a morphisms from $\pi_1(S)$ into the isometry group of
$\R^{2,1}$.  The linear part of $\rho$ determines a point in $\cT_S$ 
(\cite{mess}).

In those spacetimes Mess pointed out  a particular Cauchy surface that
is  obtained by grafting the hyperbolic surface corresponding 
to the linear holonomy along a  measured geodesic lamination.
In this way a measured geodesic lamination is associated to every
 MGHC flat spacetime (\cite{mess, benedetti-bonsante}).

Maximal flat globally hyperbolic spacetimes
come in pairs, with the elements of each pair sharing the same
holonomy: each pair contains one future complete and one past complete
spacetime. 

\begin{theorem} \label{tm:flat}
For each $\lambda_-,\lambda_+\in \cML_S$ which fill $S$, there is a unique
pair of MGHC flat spacetimes, $D_-, D_+$ with the same holonomy, respectively
past and future complete, such that $\lambda_-$ and $\lambda_+$ are the laminations
associated respectively with $D_-$ and $D_+$. 
\end{theorem}

The proof is in Appendix A. 

\subsection*{Acknowledgements} 

We wish to thank Caroline Series for several useful conversations and 
Gabriele Mondello for suggesting the argument of Proposition \ref{cont:prop}.
We are also deeply grateful to Steve Kerckhoff, who suggested many important
improvements to a previous version of the text and whose help was crucial in
making the arguments more widely understandable.

%% file: fb1.tex
\section{Preliminaries} \label{sc:2}

\subsection{Earthquakes}

According to Thurston's Earthquake Theorem \cite{kerckhoff,mess}, 
given two elements $u,v\in\mathcal T_S$, there is a unique 
$(\lambda,\mu)\in\mathcal M\mathcal L_S\times\mathcal M\mathcal L_S$ 
such that
\[
    E_l^\lambda(u)=E_r^\mu(u)=v~. 
\]
So we can consider the map
\begin{equation}\label{phi:eq}
\Phi':\mathcal T_S\times\mathcal T_S\rightarrow\mathcal M\mathcal L_S^2
\end{equation}
associating to $(u,v)$ the pair of measured lamination $(\lambda,\mu)$.
It will be clear below that this map $\Phi'$ differs from the map
$\Phi$ introduced in the previous section only by a factor $2$.

Recall that $(E^\lambda_r)^{-1}=E^\lambda_l$. 
So if $u$ is a fixed point of $E_l^\mu\circ E_l^\lambda$ then
\[
    E_l^\lambda(u)= E_r^\mu(u)
\]
which, in turn, is equivalent to
\[
   \Phi'(u,v)=(\lambda,\mu)
\]
where we have put $v= E_l^\lambda(u)$.

Conversely, if $\Phi'(u,v)=(\lambda,\mu)$ then
\[
   v=E_l^\lambda(u)=E_r^\mu(u)
\]
so that $u=E_l^\mu\circ E_l^\lambda(u)$.   
Therefore,  $E_l^\mu\circ E_l^\lambda$ admits a fixed point if and only if 
$(\lambda,\mu)$ is contained in the image of $\Phi'$. Moreover fixed 
points of $E_l^\mu\circ E_l^\lambda$ are obtained by projecting 
$\Phi'^{-1}(\lambda,\mu)$ on the first factor.
As a consequence,  
our fixed point problem can be reduced to studying the image of $\Phi'$.

\subsection{Some AdS geometry}

We briefly recall here some basic geometric properties of the anti-de Sitter
space. More details can be found e.g. in \cite{mess,mess-notes,maximal}.
 
\subsubsection*{Definition.}

The $n$-dimensional anti-de Sitter space can be defined as a quadric
in the $(n+1)$-dimensional flat space $\R^{n-1,2}$ of signature $(n-1,2)$:
$$ AdS_n = \{ x\in \R^{n-1,2}~|~ \langle x,x\rangle = -1\}~, $$
with the induced metric.
It is a geodesically complete Lorentzian manifold of constant curvature 
$-1$, with isometry group $O(n-1,2)$. It can be considered as the
Lorentzian analog of the $n$-dimensional hyperbolic space. 

One difference however is that $AdS_n$ is not simply connected,
its fundamental group is isomorphic to $\Z$, and $AdS_n$ is
homeomorphic to the product of a $(n-1)$-dimensional ball and 
a circle. It is often useful
to consider its universal cover, denoted by $\tilde{AdS}_n$ here.

\subsubsection*{Space-like, time-like and light-like vectors.}

We have mentioned that the metric induced by $\R^{n-1,2}$ is Lorentzian.
So a non-zero vector $v$ tangent to $AdS_n$ can be of three types:
\begin{itemize}
\item {\it Space-like} if $\langle v,v\rangle>0$. 
\item {\it Time-like} if $\langle v,v\rangle<0$.
\item {\it Light-like} or {\it isotropic} if $\langle v,v\rangle=0$.
\end{itemize}
This terminology originates from relativity theory, see e.g. \cite{HE}.

\subsubsection*{Geodesics and hyperplanes.}

By an elementary symmetry argument, the intersection of $AdS_n$ with 
a 2-planes $P\subset \R^{n-1,2}$ containing the origin is a geodesic $g$.
When the restriction to $P$ of the metric of $\R^{n-1,2}$ is negative
definite, the restriction to $g$ of the induced metric of $AdS_n$ is negative definite. 
Those geodesics are called time-like, and their non-zero tangent vectors are time-like. 
Conversely, time-like geodesics in $AdS_n$ are exactly 
the intersections of $AdS_n$ with the 2-dimensional planes
containing $0$ with negative definite induced metrics. It
follows that those time-like geodesics are closed, of length
equal to $2\pi$. Each of those time-like geodesics is a retract
by deformation of $AdS_n$. 

Let $P$ be a 2-dimensional plane containing $0$ in $\R^{n-1,2}$ 
on which the induced metric is of signature $(1,1)$. The intersection
of $P$ with $AdS_n$ is a complete geodesic, on which the induced
metric is positive definite -- those geodesics are called space-like.

The intersections of $AdS_n$ with the hyperplanes in $\R^{n-1,2}$
containing $0$ and of signature $(n-1,1)$ are totally geodesic,
and the induced metric is Riemannian. They are isometric to 
the hyperbolic $(n-1)$-dimensional space. Those hyperplanes are
called {\it space-like}, and their non-zero tangent vectors are
all space-like. 

The intersections of $AdS_n$ with the hyperplanes containing $0$ 
and of signature $(n-2,2)$ are also totally geodesic, but of
Lorentzian signature. They are isometric to $AdS_{n-1}$, and
are called {\it time-like} hyperplanes. Note however that their
non-zero tangent vectors are not all time-like, but they can be
either space-like, time-like or light-like (isotropic). Actually
a totally geodesic hyperplane is time-like if and only if it
contains at least one time-like tangent vector. 

\subsubsection*{Causal structure.}

We consider in the sequel a time orientation in $AdS_n$. In
other terms, we choose one of the two connected components
of the space of time-like vectors in $AdS_n$, and consider
those vectors to be future-oriented, while those in the
other connected component are past-oriented. 

We can then define the future of a subset $\Omega\subset AdS_n$, as the subset
of points in $AdS_n$ which can be obtained as the endpoint
of a future-oriented time-like curve starting from $\Omega$.
We will denote the future of $\Omega$ by $I^+(\Omega)$. 
Analogously we define the past of $\Omega$ as the set of endpoints
of past-oriented time-like curves originating from $\Omega$,
and denote it by $I^-(\Omega)$.

This notion is not too helpful in $AdS_n$, because
there are closed time-like curves.
For instance it is not
difficult to check that the future of a totally geodesic
space-like hyperplane is the whole space. 
However it is more interesting in the universal cover $\tilde{AdS}_n$.
For instance, the future of a space-like totally geodesic hyperplane
$P_0$ is a connected component of $\tilde{AdS}_n\setminus P_0$. 

\subsubsection*{Projective model, boundary at infinity}

There is a natural action of $\Z/2\Z$ on $AdS_n$, given by $x\mapsto -x$
in the quadric model above. 
The quotient space $AdS_n/(\Z/2\Z)$ has a projective model, which is often useful to
obtain a heuristic idea of its geometry. It is obtained in the same
manner as the Klein model of hyperbolic space, by projecting from the
quadric in $\R^{n-1,2}$ to the tangent hyperplane $P_0$ at one of its
points $x_0$, in the direction of $0$. This projection map, restricted
to the points $x\in AdS_n$ such that $\langle x,x_0\rangle<0$, is
projective --- it sends geodesics to line segments --- and its image
is the interior of a quadric of signature $(n-2,1)$ in $\R^n$.
It is sometimes more convenient to consider this model in the projective
space $\R P^n$, rather than in $\R^n$. One gets in this manner a projective
model of the quotient of $AdS_n$ by the ``antipodal'' action of $\Z/2\Z$.

The {\it boundary at infinity} of $AdS_n/(\Z/2\Z)$ can be defined in this
projective model, as the quadric bounding the projective model of
$AdS_n/(\Z/2\Z)$ in $\R P^n$. 
Equivalently, we can define the boundary of $AdS_n$, seen
as a quadric in $\R^{n-1,2}$, as the quotient of the space of
non-zero isotropic vectors in $\R^{n-1,2}$ by $\R_{>0}$. 

This boundary projects to $\R P^n$, and we obtain in this way
the boundary of the projective model in $\R P^n$ of the  
quotient $AdS_n/(Z/2\Z)$.

\subsubsection*{Normal vectors to hyperplanes.}

The metric on $AdS_n$ is non-degenerate. It follows that, given any
totally geodesic hyperplane $P\subset AdS_n$ and any point $x\in P$,
the orthogonal to $T_xP$ in $T_xAdS_n$ is a well-defined line in 
$T_xAdS_n$. If $P$ is oriented and non-isotropic, then we can define
its unit oriented normal vector at $x$, which is a vector $n$ such
that $\langle n,n\rangle=1$ if $P$ is time-like, and such that 
$\langle n,n\rangle =-1$ if $P$ is space-like. 

It is useful to remark that the unit normal vector field defined 
in this way is actually a parallel vector field along $P$, as can
be checked by using the quadric model of $AdS_n$ as defined above.

\subsubsection*{Angles between space-like hyperplanes.}

Let $P_1, P_2$ be two space-like totally geodesic hyperplanes in $AdS_n$,
which intersect along a codimension $2$ plane.
Then $P_1$ and $P_2$ are the intersections with $AdS_n$ of 
two hyperplanes $H_1$ and $H_2$ in $\R^{n-1,2}$ containing $0$,
each of signature $(n-1,1)$. 

We can define the angle $\theta$ between $P_1$ and $P_2$ as
the angle between $H_1$ and $H_2$ in $\R^{n-1,2}$. Let 
$N_1, N_2$ be the two unit orthogonal vectors to
$H_1$ and $H_2$ which are in the same connected component of 
$\R^{n-1,2}\setminus H_1$, then $\theta$ is the non-negative number
defined by the equation
\begin{equation}
  \label{eq:angle}
  \cosh(\theta)=|\langle N_1,N_2\rangle|~. 
\end{equation}

If $x\in P_1\cap P_2$ and $n_1,n_2$ are the unit future-oriented normals
to $P_1$ and $P_2$ at $x$, then $\theta$ can also be defined by the
fact that $\cosh(\theta)=|\langle n_1, n_2\rangle|$. 

\subsubsection*{Orthogonality between space-like and time-like hyperplanes.}

There is a well-defined notion of angle between a space-like and time-like
hyperplane in $AdS_n$ (see e.g. \cite{cpt}), but we will not really need
this notion here. What we do need, however, is the notion of orthogonality
between a time-like hyperplane and a space-like hyperplane, and also
between a time-like hyperplane and a space-like geodesic line or between
a space-like hyperplane and a time-like line.

Those notions can be defined as follows. Given a time-like hyperplane $T$
and a space-like hyperplane $S$, we say that they are orthogonal at a point
$x\in P\cap S$ if, at $x$, the unit vector normal to $T$ is orthogonal to
the unit vector normal to $S$. Since the normal vector fields are parallel
along $T$ and along $S$, it then follows that $T$ and $S$ are also
orthogonal at any other intersection point, and we will say simply that
they are orthogonal. 

Given now a space-like geodesic line $D$ intersecting $T$, we say that
$D$ is orthogonal to $T$ if, at the intersection point, $D$ is parallel to 
the unit vector normal to $T$. 
In the same way, a time-like line $D$ is orthogonal to a space-like
hyperplane $S$ if, at the intersection point, $D$ is parallel to the
unit vector normal to $S$.

\subsubsection*{AdS vs de Sitter.} 

It can be useful to recall that the de Sitter $n$-dimensional space
$dS_n$ is defined as a quadric in the Minkowski $(n+1)$-dimensional
space: 
$$ dS_n = \{ x\in \R^{n,1}~|~ \langle x,x\rangle =1\}~. $$
It is also a Lorentzian constant curvature space, but of curvature 
$1$ rather than $-1$. The de Sitter space is dual to the hyperbolic
space, rather than analogous to it, see e.g. \cite{RH}. 

For $n=2$, however, $AdS_2$ and $dS_2$ are very similar, since 
$AdS_2$ is isometric to $dS_2$ with a reversed sign (in other
terms, there is a map from $AdS_2$ to $dS_2$ which only changes
the sign of the scalar product). 

\subsubsection*{$AdS_3$.}

The $3$-dimensional AdS space, $AdS_3$, has some very specific
properties, see \cite{mess}. One reason for this is that $AdS_3/(\Z/2\Z)$ is
isometric to $PSL_2(\R)$ with its Killing metric. As a consequence,
$PSL_2(\R)$ acts isometrically on $AdS_3/(\Z/2\Z)$ by left and right multiplication.
This identifies $PSL_2(\R)\times PSL_2(\R)$ with the identity component
of the isometry group of $AdS_3/(\Z/2\Z)$. 

\subsection{AdS geometry and hyperbolic surfaces}

G. Mess, in his celebrated work \cite{mess}, showed that  the map $\Phi'$ 
is meaningful in AdS context. 
We recall here how Mess connected the map $\Phi'$ to AdS geometry. 
This connection will play a fundamental role in the rest of the paper.

The relation between AdS and hyperbolic surfaces arises from the fact 
that AdS manifolds are locally modeled on the model $AdS_3/(\Z/2\Z)$ 
whose isometry group is $PSL_2(\mathbb R)\times PSL_2(\mathbb R)$.  
Moreover, the boundary at infinity of $AdS_3/(\Z/2\Z)$ is naturally identified to 
$\partial\mathbb H^2\times\partial\mathbb H^2$. The action of 
$PSL_2(\mathbb R)\times PSL_2(\mathbb R)$ on $AdS_3/(\Z/2\Z)$ extends to the boundary. In fact,
the action  on the boundary coincides with the product action.

Given two Fuchsian representations $\rho_l,\rho_r:\pi_1(S)\rightarrow
PSL_2(\mathbb R)$, there exists a unique homeomorphism $\phi$ of
$\partial \mathbb H^2$ conjugating their actions on $\R P^1=\partial
\mathbb H^2$ (\cite{ahlfors:riemann}).

Notice that the pair $\rho=(\rho_l,\rho_r)$ can be regarded as a
representation in $PSL_2(\mathbb R)\times PSL_2(\mathbb R)$. The graph $\Gamma_\rho$ of
$\phi$ turns out to be an invariant subset of $\partial\mathbb
H^2\times\partial\mathbb H^2$.
In fact, it is the minimal non-empty invariant closed subset, like the
limit set of a quasifuchsian representation.

Mess showed that there exists an open convex subset
$\Omega=\Omega(\rho)$ in $AdS_3/(\Z/2\Z)$ that is invariant under $\rho$ such
that the action of $\rho$ on $\Omega$ is free and properly
discontinuous and the quotient $M(\rho)=\Omega/\rho$ is a globally
hyperbolic AdS-manifold diffeomorphic to $S\times\mathbb R$.  In fact,
such $\Omega$ can be chosen maximally, in the sense that any other
domain satisfying the same properties is contained in $\Omega$.

Spacetimes $M(\rho)$ constructed in this way completely classify globally
hyperbolic AdS structures on $S\times\mathbb R$.

\begin{prop}\cite{mess}
Let $M$ be a globally hyperbolic AdS spacetime diffeomorphic to 
$S\times\mathbb R$. Then there is a pair of Fuchsian representations
$\rho=(\rho_l,\rho_r)$ of $\pi_1(S)$ such that
$M$ isometrically embeds into $M(\rho)$.
\end{prop}   

Given a pair of Fuchsian representations $\rho=(\rho_l,\rho_r)$, the
set $\Gamma_\rho$ is contained in the closure of $\Omega(\rho)$ in
$\overline{AdS_3/(\Z/2\Z)}$, see \cite{mess}. 
Thus the convex hull $C(\tilde M(\rho))$ of
$\Gamma_\rho$ in $AdS_3/(\Z/2\Z)$ is contained in $\Omega$ and it projects to a
convex subset $C(M(\rho))$ of $M(\rho)$.  This subset turns out to be
the \emph{convex core} of $M(\rho)$, in the sense that it is a
non-empty convex closed subset  of $M(\rho)$ that is
contained in every non-empty closed convex subset of $M(\rho)$. Here
we say that $K\subset M(\rho)$ is {\it convex} if any any geodesic
segment with endpoints in $K$ is contained in $K$.  With this
definition any non-empty convex subset is a deformation retract of
$M(\rho)$.

If $\rho_l$ and $\rho_r$ are conjugate, then $C(M(\rho))$ is a totally
geodesic
space-like surface homeomorphic to $S$, and the restriction of the
metric on $C(M(\rho))$ makes it isometric to $\mathbb H^2/\rho_l$.  If
$\rho_l$ and $\rho_r$ represent two different points in $\mathcal T$
then $C(M(\rho))$ is homeomorphic to $S\times [-1,1]$ and it is a
strong deformation retract of $M(\rho)$.  The boundary components of
$C(M(\rho))$ are achronal $\mathrm C^{0,1}$-surfaces (achronal means
that time-like paths meet these sets in at most one point).  We denote
by $\dr_-C(M(\rho))$ (resp. $\dr_+C(M(\rho))$) the past (resp. future)
component of $\partial C(M(\rho))$.

The minimality condition implies that $\dr_-C(M(\rho))$ and $\dr_+C(M(\rho))$
are totally geodesic space-like surfaces bent along geodesic laminations $L_-$
and $L_+$ respectively. Since  space-like planes in $AdS_3$ and in  $AdS_3/(\Z/2\Z)$ are
isometric to $\mathbb H^2$, the boundary components of
$C(M(\rho))$ are equipped with a hyperbolic structure. Moreover, it is
possible to equip $L_-$ and $L_+$ with transverse measures $\mu_-$ and
$\mu_+$ that measure the amount of bending \cite{mess}.
 
Mess \cite{mess} noticed that the right earthquake along
$2\lambda_+=(L_+,2\mu_+)$ transforms $\mathbb H^2/\rho_l$ into
$\mathbb H^2/\rho_r$ and analogously
the left earthquake along $2\lambda_-=(L_-,2\mu_-)$ transforms $\mathbb
H^2/\rho_l$ into $\mathbb H^2/\rho_r$.
 With our notation, we have
 \begin{equation}\label{mm:eq}
   \Phi'(\rho_l,\rho_r)=(2\lambda_-,2\lambda_+)
\end{equation}
so that $\Phi'$ appears as twice the map $\Phi$ in Section 1.

Thus, the problem of determining the image of $\Phi'$ is equivalent to
the problem of determining which pairs of laminations can be realized
as bending laminations of some AdS globally hyperbolic structure on
$S\times\mathbb R$.  A first constraint for $(\lambda_-,\lambda_+)$ to
lie in the image of $\Phi'$ is that they have to fill the surface.

\begin{lemma}
If $u\neq v$ then $\Phi'(u,v)$ is a pair of measured laminations that fill the surface $S$.
\end{lemma}
\begin{proof}

By Formula (\ref{mm:eq}), the statement is equivalent to the fact that the bending
laminations of every GHMC AdS spacetime homeomorphic to $S\times \mathbb R$ 
fill $S$.

This fact can be proved along the same line as the analogous result for quasifuchsian manifold.
Suppose that a curve $c$  meets neither $\lambda_+$ nor $\lambda_-$. 
Then $c$ can be realized as an unbroken geodesic both in $\partial_+C(M(\rho))$ and in 
$\partial_-C(M(\rho))$. 
Thus two geodesic representatives of $c$ should exists in $M$. But this contradicts the fact 
that for any free homotopy class of loop in $M$ there is a unique geodesic representative (see 
\cite{mess,benedetti-bonsante}). 
\end{proof}

\begin{remark}
The description of globally hyperbolic AdS manifolds, including the
geometry of the boundary of the convex core, is strongly reminiscent
of quasifuchsian hyperbolic manifolds.  In the quasifuchsian case the
convex core also exists; its boundary has an induced metric which is
hyperbolic and it is pleated along a measured lamination. However the
measured bending laminations on the boundary satisfy an additional
condition: the weight of any closed curve is less than $\pi$. 

The reason why such a restriction is not needed in AdS case lies on
the fact that in AdS geometry, and more generally in Lorentzian
geometry, the angle between two space-like totally geodesic planes is a
well defined number in $[0,+\infty)$ (see Equation (\ref{eq:angle})).
\end{remark}

\section{Existence of fixed points for small laminations} \label{sc:3}

\subsection{General setup}

Let $\Delta$ denote the diagonal of $\mathcal T\times\mathcal T$.  In
this section we will prove that there is a neighborhood $\hat U$ of
$\Delta$ in $\mathcal T\times\mathcal T$ such that the restriction of
$\Phi$ to $U=\hat U\setminus\Delta$ is a homeomorphism onto an open
set $V\subset\mathcal F\mathcal M\mathcal L_S$.  Moreover, we will show
that $\hat U$ can be chosen such that for any $(\lambda,\mu)\in \cFML_S$
we have that $(t\lambda, t\mu)\in V$ for $t$ sufficiently small.

\subsection{Infinitesimal earthquakes and the length function}

For a given measured geodesic lamination $\lambda$ on $S$, the
following semi-group law holds:
\[
     E_l^{t\lambda} \circ E_l^{s\lambda}= E_l^{(t+s)\lambda}\,.
\]
We consider the \emph{infinitesimal left earthquake}
$e_l^\lambda(u)=\frac{\mathrm d\,}{\mathrm dt}E_l^{t\lambda}(u)$, which
is a vector field on $\mathcal T$.  It follows from the semi-group relation
that $(E_l^{t\lambda})_{t\geq 0}$ corresponds to the flow of
$e_l^\lambda$ for positive time.

Analogously one can define the \emph{infinitesimal right earthquake}
$e_r^\lambda$ and remark that $(E_r^{t\lambda})_{t\geq 0}$ is the flow
for positive time of this field. Moreover, since
$E_r^{\lambda}=(E_l^{\lambda})^{-1}$ we have that
$e_r^\lambda=-e_l^\lambda$ for every $\lambda$.

A basic standard property is that these fields continuously depend on
the lamination $\lambda$.  For the convenience of the reader, we will
sketch the proof of this fact.

\begin{prop}\label{cont:prop}
Let $\Gamma(T\mathcal T)$ denote the space of vector fields on $\mathcal T$
equipped with the $\mathrm C^\infty$-topology.
The map
\[
    \mathcal M\mathcal L_S\ni\lambda\mapsto e_l^\lambda\in\Gamma(T\mathcal T)
\]
is continuous.
\end{prop}
\begin{proof}
Let $l_\gamma$ denote the length function of the lamination $\lambda$, which is
a smooth function on $\mathcal T$.
A celebrated result of Wolpert states that the field
$e_l^\lambda$ is the symplectic gradient of  $l_\lambda$ with respect
to the Weil-Petersson form (see \cite{wolpert-formula}).  
Thus it is sufficient to prove that the map
\[
    \mathcal M\mathcal L_S\ni\lambda\mapsto l_\lambda\in\mathrm
    C^\infty(\mathcal T)
 \]
is continuous.
  
Consider $\mathcal T$ as the Fuchsian locus of the space of
quasifuchsian metrics $\mathcal Q\mathcal F$ on $S\times\mathbb R$,
that is in a natural way a complex manifold.
Bonahon \cite{bonahon-toulouse} proved that the function $l_\lambda$
can be extended to a holomorphic function --- still denoted by
$l_\lambda$ --- on $\mathcal Q\mathcal F$.
 
Let $\mathcal O(\mathcal Q\mathcal F)$ be the space of holomorphic functions on
$\mathcal Q\mathcal F$.
The map
 \[
   l: \mathcal M\mathcal L_S\ni\lambda\mapsto l_\lambda\in\mathcal
   O(\mathcal Q\mathcal F)
 \]
 is locally bounded and continuous with respect to the pointwise
 topology on $\mathcal O(\mathcal Q\mathcal F)$. Montel's theorem
 shows that it is continuous if $\mathcal O(\mathcal Q\mathcal F)$ is
 equipped with the $\mathrm C^\infty$-topology.
\end{proof}

The results of this section will be based on a transversality argument
that has been developed
in \cite{kerckhoff2, series-minima, bonahon-almost}. 

\begin{prop}\label{kerck:prop} \cite{bonahon-almost}
Let $\lambda,\mu\in\mathcal M\mathcal L_S$ be two measured
laminations. The intersection between $e_l^\lambda$ and $e_r^\mu$,
considered as submanifolds in the total space of the bundle $T\cT$, is
transverse.  Moreover if $\lambda$ and $\mu$ fill up the surface then
these sections meet in exactly one point $k_0(\lambda,\mu)$. Otherwise
they are disjoint.
\end{prop}

\begin{remark}
The point $k_0(\lambda,\mu)$ continuously depends on $\lambda$ and
$\mu$ and it is the unique minimum point for the proper function
$l_\lambda+l_\mu$. The fact that $e_l^\lambda$ and $e_r^\mu$
transversely intersect in $k_0(\lambda,\mu)$ is a consequence of the
fact that the Hessian of $l_\lambda+l_\mu$ at $k_0(\lambda,\mu)$ is
positive definite (see \cite{kerckhoff2}).
\end{remark}

\subsection{Bending laminations near the Fuchsian locus}

Through this section, for $(\lambda,\mu)\in\cFML_S$, we denote by
$k_0(\lambda,\mu)$ the unique point of $\mathcal T$ 
where $e^\lambda_l$ is equal to $e^\mu_r$.

Using the transversality result of Proposition \ref{kerck:prop}, 
we prove that for every pair of laminations $\lambda,\mu$ that fill
the surface the composition $E_l^{t\mu}\circ E_l^{t\lambda}$
admits a fixed point for small $t$. 

\begin{prop}\label{loc:prop}
Let  $\lambda_0$ and $\mu_0$ be laminations that fill up the surface
and $u_0=k_0(\lambda_0,\mu_0)$. 
There exist  $\epsilon>0$, a neighborhood $V$ of
$(\lambda_0,\mu_0)\in\mathcal F\mathcal M\mathcal L_S$, a
neighborhood $U$ of $u_0$, and a continuous map
\[
   k:V\times[0,\epsilon)\rightarrow U
\]
such that:
\begin{itemize}
\item For $t>0$, $k(\lambda,\mu,t)$ is the unique fixed point
of $E^{t\mu}_l\circ E^{t\lambda}_l$ lying in $U$. 
\item For all $(\lambda,\mu)\in V$, $k(\lambda,\mu,0)=k_0(\lambda,\mu)$. 
\end{itemize}
\end{prop}

The first idea to prove this proposition could be to consider 
the function $E^\lambda_l\circ E^\mu_l:\cT_S\rightarrow \cT_S$, 
considered as depending on the parameter $(\lambda,\mu)\in \cML_S\times \cML_S$,
prove that it has a fixed point as $(\lambda,\mu)\rightarrow 0$, and apply
the Implicit Function Theorem. This is broadly speaking the argument we use,
but some care is needed. One reason  is that the map is 
critical as $(\lambda,\mu)\rightarrow 0$, so that a blow-up procedure is
necessary. Another, related reason is that $(0,0)\not\in \cFML_S$.

A more exact description of the argument of the proof is as follows.
First, for any $(\lambda,\mu)\in\mathcal F\mathcal M\mathcal L_S$ and $t\geq 0$
we point out a smooth function 
$\phi_{\lambda,\mu,t}:\mathcal T\rightarrow\mathbb R^{6g-6}$ such that
\begin{itemize}
\item $\phi_{\lambda,\mu,0}^{-1}(0)$  
contains only $k_0(\lambda,\mu)$.
\item $\phi_{\lambda,\mu,t}^{-1}(0)$  
is the set of fixed points of $E^{t\mu}_l\circ E^{t\lambda}_l$ when $t>0$.
\end{itemize}

In particular, $u_0$ is a solution of the equation
$\phi_{\lambda_0,\mu_0,0}(u)=0$.  Using the transversality given by
Proposition~\ref{kerck:prop} we show that the differential
$\big(d\phi_{\lambda_0,\mu_0,0}\big)_{u_0}$ is not singular.  The
proof is then concluded by a simple application of the Implict
Function Theorem.  We state for the reader's convenience the Implicit
Function Theorem in the form we will use below.

\begin{lemma}\label{impl:lem} 
Let $X$ be a topological space and $M$ be a differentiable manifold of
dimension $n$.  Consider a family of smooth maps
$\{\phi_x:M\rightarrow\mathbb R^n\}_{x\in X}$ indexed by the elements
of $X$ so that the induced map
\[
     X\ni x\mapsto \phi_x\in\mathrm C^{\infty}(M,\mathbb R^n)
\]  
is continuous.

Suppose that for some $x_0\in X$ and $u_0\in \mathbb R^n$
we have $\phi_{x_0}(u_0)=0$ and $d(\phi_{x_0})_{u_0}$ is not singular.
Then, there is a neighborhood $V$ of $x_0$ and a neighborhood $U$ of $u_0$
such that for any $x\in V$ the equation
\[
     \phi_x(u)=0
\]
admits a unique solution $u(x)$ lying in $U$.
Moreover the map 
\[
   V\ni x\mapsto u(x)\in U
\]
is continuous.
\end{lemma}

\begin{proof}[Proof of Proposition \ref{loc:prop}]
   
Finding a fixed point of $E_l^{t\lambda}\circ E_l^{t\mu}$ is
equivalent to finding a point $u$ such that
\begin{equation}\label{fx:eq}
    E_l^{t\lambda}(u)=E_r^{t\mu} (u)~.
\end{equation}
Let us fix a global diffeomorphism
\[
  \mathbf{x}:\mathcal T\rightarrow\mathbb R^{6g-6}\,.
\]
(We could for instance use the Fenchel-Nielsen coordinate for some pant
decomposition of $S$, but the choice is irrelevant here.)

Consider the $\mathbb R^{6g-6}$-valued functions on $\mathcal T$
defined by  $\mathcal E_l^{t\lambda}=\mathbf{x}\circ
E_l^{t\lambda}$ and $\mathcal E_r^{t\mu}=\mathbf{x}\circ
E_r^{t\lambda}$.


Notice that $\mathcal E^{t\lambda}_l-\mathcal E^{t\mu}_r$ forms
a family of smooth maps from $\mathcal T$ to $\mathbb R^{6g-6}$
indexed by $\mathcal F\mathcal M\mathcal L_S\times \mathcal[0,\infty)$.
Unfortunately, for $t=0$, $\mathcal
E_l^{t\lambda}=\mathcal E_r^{t\mu}$, so the difference is $0$
everywhere and we cannot apply Lemma \ref{impl:lem}.

For this reason we blow up the difference.  More precisely, we
consider the family of maps
\[
   \phi_{\lambda,\mu, t}(u)= \frac{\mathcal E_l^{t\lambda}(u)-\mathcal
     E_r^{t\mu}(u)}{t}~.
\]
Clearly, if we fix $(\lambda,\mu)$ and $u\in\mathcal T$ then $
\lim_{t\rightarrow
  0}\phi_{\lambda,\mu,t}(u)=d\mathbf{x}_u(e^\lambda_l(u)-e^\mu_r(u))$,
where we are implicitly using the canonical identification of
$T_{\mathbf{x}(u)}\mathbb R^{6g-6}$ with $\mathbb R^{6g-6}$.

Moreover, putting 
\begin{equation}\label{ext:eq}
\phi_{\lambda,\mu,0}(u)=d\mathbf{x}_u(e^\lambda_l(u)-e^\mu_r(u))\,,
\end{equation}
Proposition \ref{cont:prop} implies that the induced application
\[
  \phi:\mathcal F\mathcal M\mathcal L_S\times [0,+\infty)\rightarrow
    C^\infty(\cT,\mathbb R^{6g-6})
\]
is continuous.

By the choice of $u_0$ we have $\phi_{\lambda,\mu,0}(u_0)=0$.  On the
other hand, since $d\mathbf x$ provides a trivialization of the
tangent bundle of $\mathcal T$, Proposition \ref{kerck:prop} implies
that the graphs of the functions
$$X_0(u)=d\mathbf x_u(e^{\lambda_0}_l)~,\qquad Y_0(u)=d\mathbf x_u(e^{\mu_0}_r)$$ 
transversely meet  over $u_0$. Since $\phi_{\lambda_0,\mu_0,0}=X_0-Y_0$ we  
conclude that
$d\big(\phi_{\lambda_0,\mu_0,0}\big)_{u_0}$ is non-singular.

From Lemma~\ref{impl:lem} there is a neighborhood $V$ of
$(\lambda_0,\mu_0)$ in $\cFML_S$, $\epsilon>0$, a neighborhood $U$ of $u_0$,
and a continuous function
\[
 k:V\times[0,\epsilon)\rightarrow U
\]
such that $k(\lambda,\mu,t)$ is the unique point in $U$ such that
\[
\phi_{\lambda,\mu,t}(k(\lambda,\mu,t))=0\,.
\]

Therefore, for $t>0$, $k(\lambda,\mu, t)$ is the unique fixed point for
$E^{t\mu}_l\circ E^{t\lambda}_l$ lying in $U$.
On the other hand,  Equation (\ref{ext:eq}) implies that
$k(\lambda,\mu,0)=k_0(\lambda,\mu)$.
\end{proof}

Let $\cFML_1$ be a hypersurface in $\cFML_S$ that intersects every ray
of $\cFML_S$ exactly once.  
By Proposition \ref{loc:prop}, there exists an open
covering $\{V_i\}_{i\in I}$ of $\cFML_1$ and a family of maps
\begin{equation}\label{loc:eq}
     k^{(i)}:V_i\times[0,\epsilon_i)\rightarrow U_i\subset\mathcal T
\end{equation}
such that, for all $i\in I$,
\begin{enumerate}
\item $k^{(i)}(\lambda,\mu,0)=k_0(\lambda,\mu)$.
\item If $t>0$, then $k^{(i)}(\lambda,\mu,t)$ is the unique fixed point
of $E^{t\mu}_l\circ E^{t\lambda}_l$ lying in $U_i$.
\end{enumerate}

The following lemma shows how these maps can be glued to a global map.

\begin{lemma}\label{glob:lemma}
There are an open neighborhood $V^*$ of $\cFML_1\times\{0\}$ in
$\cFML_1\times[0,+\infty)$  and a continuous map 
$k^*:V^*\rightarrow\mathcal T$  such that:
\begin{enumerate}
\item
If $(\lambda,\mu, t_0)\in V^*$ then $(\lambda,\mu, t)\in V^*$ for
every $t\in[0,t_0]$.
\item
$k^*(\lambda,\mu,0)=k_0(\lambda,\mu)$.
\item For $t>0$
$k^*(\lambda,\mu,t)$ is a fixed
point  for $E^{t\mu}_l\circ E^{t\lambda}_l$.
\end{enumerate}
\end{lemma}

\begin{proof}
Let $k^{(i)}:V_i\times[0,\epsilon_i)\rightarrow U_i\subset \cT$ be the
family of maps defined in (\ref{loc:eq}) indexed by $i\in I$.  
The uniqueness property of the function $k^{(\cdot)}$ shows that, for $i$ and $j$ in
$I$, the subset of $\big(V_i\times[0,\epsilon_i)\big)\cap \big(V_j\times[0,\epsilon_j)\big)$ 
where $k^{(i)}$ and $k^{(j)}$ concide is open. 
Since $k^{(i)}=k^{(j)}$ on $V_i\cap V_j \times\{0\}$ the
corresponding maps $k^{(i)}$ and $k^{(j)}$ coincide on  the whole
$\big(
V_i\times[0,\epsilon_i)\big)\cap\big(V_j\cap[0,\epsilon_j)\big)$.
It follows that on the set
\[
     V^*=\bigcup_i V_i\times[0,\epsilon_i)
\]
the map $k^*$ can be defined by gluing the maps $k^{(i)}$.
\end{proof}

Let 
\[
   \pi:\cFML_1\times(0,+\infty)\rightarrow\cFML_S
\]
be the homeomorphism defined by $\pi(\lambda,\mu,t)=(t\lambda, t\mu)$.
Let us  define $V=\pi\big(V^*\setminus(\cFML_1\times\{0\})\big)$
where $V^*$ is the open subset defined in Lemma \ref{glob:lemma}.
On the set $V$ we can consider the map  $k=k^*\circ \pi^{-1}$.
Notice that $k(\lambda,\mu)$ is a fixed point of $E^{\mu}_l\circ E^{\lambda}_l$.
By means of this map, we can construct on $V$ a right inverse of 
the map $\Phi'$ defined in (\ref{phi:eq}).

\begin{cor} \label{cr:loc}
The open set $V$ verifies the following properties.
\begin{enumerate}
\item
For every $(\lambda,\mu)\in\cFML_S$ there is $t>0$ such that
$(t\lambda,t\mu)\in V$.
 \item 
If $(\lambda,\mu)\in V$ then $(t\lambda,t\mu)\in V$ for every $t\in(0,1)$.
  \item
There is a continuous map $\sigma:V\rightarrow \mathcal T\times\mathcal T$ 
that is a right inverse for the map $\Phi'$. Moreover
\[
   \lim_{t\rightarrow 0}\sigma(t\lambda, t\mu)=(k_0(\lambda,\mu),
   k_0(\lambda,\mu))\,.
\]

\end{enumerate}
\end{cor}

\begin{proof} 
The first two properties follow directly from the definition on $V$.
Using the fact that $k$ sends $(\lambda,\mu)$ to a fixed point of
$E^\mu_l\circ E^\lambda_l$, the map $\sigma$ can be defined by putting
\begin{equation}
    \sigma(\lambda,\mu)=(k(\lambda,\mu), E^\lambda_l(k(\lambda,\mu)))\,.
\end{equation}
\end{proof}

\subsection{All metrics near the Fuchsian locus are obtained} 

Through this section we use the same notations as in the previous
section.  In particular we consider the map $\sigma$ constructed in
Corollary \ref{cr:loc}.  This map is clearly injective, so by the
Theorem of the Invariance of Domain the image of $\sigma$ is an open
set $U$ in $\mathcal T\times\mathcal T$.
Notice that the restriction map $\Phi'|_U$ is a homeomorphism of $U$
onto $V$.

\begin{prop}\label{pr:U}
$\hat U=U\cup\Delta$ is an open neighborhood of the diagonal $\Delta$
  in $\mathcal T\times\mathcal T$.
\end{prop}

To prove Proposition \ref{pr:U}, we need two technical results about
the behavior of the map $\Phi'$ near the diagonal.

\begin{lemma}\label{tec1:lem}
Let $(u_k,v_k)_{k\in \N}\in\mathcal T^2$ be a sequence converging to
$(u,u)\in\Delta\subset\mathcal T^2$ and $(t_k)_{k\in \N}$ be a
sequence of positive numbers such that, putting $\Phi'(u_k, v_k)=(t_k\lambda_k,
t_k\mu_k)$, we have that $(\lambda_k)$ converges to a measured lamination $\lambda\neq 0$.
Then the sequence of measured laminations $(\mu_k)$ also converges to a measured 
lamination $\mu$. Moreover
$\lambda$ and $\mu$ fill the surface 
and $u$ is the point $k_0(\lambda,\mu)$ where $e^\lambda_l$ 
and $e^\mu_r$ meet.
\end{lemma}
\begin{proof}
We have
\[
   E^{t_k\lambda_k}_l(u_k)=E^{t_k\mu_k}_r(u_k)
\]
so by considering the Taylor expansion we deduce that
\[
   e^{\lambda_k}_l(u_k)-e^{\mu_k}_r(u_k)=O(t_k)\,.
\]
Taking the limit for $k\rightarrow+\infty$ we have that
\[
   \lim_{k\rightarrow +\infty}e^{\mu_k}_r(u_k)=e^\lambda_l(u)\,.
\]
Since the map $\mathcal M\mathcal L\times\mathcal T\ni(\mu, u)\mapsto e^\mu_r(u)\in T\mathcal T$
is a homeomorphism (see \cite{kerckhoff2})
we deduce that $(\mu_k)$ converges to a measured lamination $\mu$ such that
\[
     e^{\lambda}_l(u)=e^{\mu}_r(u)\,.
\]      
Thus, since we are assuming $\lambda\neq 0$, we have that $\lambda$
and $\mu$ fill the surface and $u=k_0(\lambda,\mu)$.
\end{proof}

\begin{lemma}\label{tec2:lem}
If $(u_k,v_k)\in\mathcal T^2$ converges to $(u,u)$, then there exists
a sequence of positive numbers $(t_k)_{k\in \N}$ such that if we
put $\Phi'(u_k, v_k)=(t_k\lambda_k, t_k\mu_k)$ then 
$\{(\lambda_k,\mu_k)\}$ is precompact in $\mathcal F\mathcal M\mathcal L_S$.
\end{lemma}

This lemma states that the two factors of $\Phi'(u_k,v_k)$ 
approach zero with the same speed.

\begin{proof}
We fix a sequence $t_k\rightarrow 0$ such that if we put
$\Phi'(u_k,v_k)=(t_k\lambda_k, t_k\mu_k)$ we have that the sequence
$\lambda_k$ is precompact in $\mathcal M\mathcal L\setminus\{0\}$.
By Lemma \ref{tec1:lem}, for every subsequence $k_j$ such that
$\lambda_{k_j}$ converges we have that $(\lambda_{k_j},\mu_{k_j})$
converges to a point in $\mathcal F\mathcal M\mathcal L_S$.
This proves that the set $\{(\lambda_k,\mu_k)\}$ is precompact in
$\mathcal F\mathcal M\mathcal L_S$.
\end{proof}

We can prove now Proposition \ref{pr:U}.

\begin{proof}[Proof of Proposition \ref{pr:U}]
Suppose by contradiction that $\hat U$ is not a neighborhood of
$\Delta$.  There is a sequence $(u_k, v_k)\in\mathcal T\times\mathcal
T$ converging to some point $(u,u)\in\Delta$ such that
$(u_k,v_k)\notin\hat U$.  

By Lemma \ref{tec2:lem}, there is a sequence of positive
numbers $t_k\rightarrow 0$ such that if we put
$\Phi'(u_k,v_k)=(t_k\lambda_k, t_k\mu_k)$ then --- up to passing to a
subsequence --- $(\lambda_k,\mu_k)\rightarrow(\lambda,\mu)\in\mathcal
F\mathcal M\mathcal L_S$. By property (1) of Corollary \ref{cr:loc},
there is $\epsilon>0$ such that $(t\lambda, t\mu)\in V$ for
$t<2\epsilon$.  Since $V$ is open we have that $(\epsilon\lambda_k,
\epsilon\mu_k)$ lies in $V$ for $k$ sufficiently large. On the other
hand, for $k$ large enough, $t_k<\epsilon$ and, by property (2) of
Corollary \ref{cr:loc} we deduce that $\Phi'(u_k,v_k)=(t_k\lambda_k,
t_k\mu_k)$ lies in $V$. Thus $(u_k,v_k)=\sigma(t_k\lambda_k, t_k\mu_k)$ and
this contradicts our assuption on the sequence $(u_k,v_k)$.
\end{proof}
 

%% file: jms2.tex
\section{An upper bound on the length of laminations} \label{sc:4}

\subsection{The main estimate.}

The goal of this section is to prove the following key estimate on 
pairs of laminations that fill a surface. 

\begin{prop} \label{pr:main}
There exist constants $\epsilon_0, h_0>0$ (depending only on the genus
of $S$) as follows. Let $(\lambda,\mu)\in \cFML_S$, and let $g\in
\cT_S$ be such that $E_r^\lambda(g) = E_l^\mu(g)$. Then
$i(\lambda,\mu)\geq \epsilon_0 l_g(\lambda) \min(l_g(\lambda), h_0)$.
\end{prop}

This proposition will be used twice below.  In particular, note that the inequality
it contains is qualitatively different, depending on whether $l_g(\lambda)$
is small or not.  When it is small the inequality is quadratic in $l_g(\lambda)$ and that
fact will be important in Section 6 where Proposition \ref{pr:main} will be a key tool for
understanding the behavior of the pleating lamination in the neighborhood of the
Fuchsian locus.  On the other hand, it is only linear in $l_g(\lambda)$ when
$l_g(\lambda)$ is large.
In that form, it will play a key role in Section 5 in the proof of the main
compactness statement, Proposition \ref{pr:proper}.  In fact, in that section we will use
the following simple consequence of Proposition \ref{pr:main}.

\begin{cor} \label{cr:obvious}
For all $C_1>0$ there exists $C_2>0$ such that, 
if $i(\lambda,\mu)\leq C_1$ and if $E_l^\lambda(g) = E_r^\mu(g)$,
then $l_g(\lambda)\leq C_2$ and $l_g(\mu)\leq C_2$.
\end{cor}

The proof of Proposition \ref{pr:main} is based on 
the geometry of 3-dimensional AdS manifolds. The fact that $E_r^\lambda(g) = E_l^\mu(g)$
means that $g$ is the left representation of a globally hyperbolic AdS manifold
for which the bending lamination on the upper boundary of the convex core is
$\lambda_+=\lambda/2$, while the bending lamination on the lower boundary of the convex core
is $\lambda_-=\mu/2$. The left representation is then $E_l^\lambda(g)$.

Proposition 4.1 will be proved by bounding $i(\lambda_{+},\lambda_-)$ in terms
of the length of $\lambda_{+}$ for the induced metric on the upper boundary of the convex core.  
This will be done by first considering the case when $\lambda_{+}$ is a weighted simple closed geodesic.  
In this case, the estimate will come from the area of a timelike geodesic annulus 
which hits the top in this closed geodesic and hits the bottom in a curve transverse to $\lambda_-$.   
In particular we want to measure the ``width'' of the convex core of $M$ in certain directions.

\subsection{The width of the convex core of a MGHC AdS spacetime}

Let $M$ be a globally hyperbolic AdS spacetime homeomorphic to
$S\times \mathbb R$ whose upper lamination $\lambda_+$ is a weighted
curve $(c,w)$.  We fix $x\in c$ and we
consider a timelike geodesic $\tau$ through $x$
that is orthogonal to a face containing $x$.
In this section we find an estimate from below for the length of
$\tau\cap C(M)$.  This estimate will be a key point in the proof of
Proposition \ref{pr:main}.

It is based on the next two technical statements in AdS geometry. Recall (see
Section 2.3) that space-like surfaces in $AdS_3$ have a hyperbolic
induced metric, as do pleated surfaces in $H^3$. However there are some
diferences between the geometry of pleated surfaces in the AdS and
the hyperbolic case, as displayed in the following statement.

\begin{sublemma}\label{sublemma}
Let $\Delta_+$ be a convex surfaces obtained by bending a spacelike plane 
in $\tilde{AdS}_3$ along a locally finite measured lamination $\tilde\lambda_+$. 
Let $\sigma$ be any geodesic path in $\Delta_+$ joining a point
$x\in \tilde\lambda_+$ to some point $y$. Let $P$ be the space-like plane through
$x$ extending the face of $\Delta_+$ that does not meet
$\sigma$. Analogously let $Q$ be the spacelike plane extending the face
containing $y$ (if $y$ lies on a bending line, we choose the face that
does not meet $\sigma$).  Then $P$ and $Q$ meet along a spacelike line $r$.
Moreover the following two properties hold:
\begin{itemize}
\item
If $\alpha(P,Q)$ denotes the angle between $P$ and $Q$ (as defined in
Section 2.2), 
\begin{equation}\label{prop1:eq}
\alpha(P,Q)\geq i(\sigma,\tilde\lambda_+)\,.
\end{equation} 
\item
If $d_P(x,r)$ is the distance from $x$ to $r$ on the plane $P$ and
$d_Q(y,r)$ is the distance from $y$ to $r$ on $Q$, then
\begin{equation}\label{prop2:eq}
d_P(x,r)+d_Q(y,r)\leq l(\sigma)~,
\end{equation}
where $l(\sigma)$ is the length of $\sigma$.
\end{itemize}
\end{sublemma}

\begin{proof}
Suppose that $P$ and $Q$ are disjoint. Then, up to exchanging $P$ and
$Q$, $P$ lies in the future of $Q$. On the other hand, $\Delta_+$ is
contained in $\overline{I^-(Q)}$ and meets $P$ at $x$ (recall that 
$I^-(Q)$ is the past of $Q$, as defined in Section 2). This gives a
contradiction.

In order to prove (\ref{prop1:eq}) first notice that we can easily
reduce to the case where $y$ lies on a bending line. Indeed let $y'$
be the point on $\sigma\cap\tilde\lambda_+$ contained in the same face as
$y$: if (\ref{prop1:eq}) holds for the segment of $\sigma$ joining $x$ to
$y'$, then it can be easily checked that it holds for $\sigma$.

Assuming that $y$ lies on some bending line, we proceed by induction on
the number of bending lines between $x$ and $y$.
Let $x_1=x,\ldots, x_n=y$ be the intersection points of $\sigma$ with
$\tilde\lambda_+$.  For the case $n=2$, we consider the plane $R$
extending the face through $x$ containing $\sigma$.  Lemma 6.15 of
\cite{benedetti-bonsante} states
$\alpha(P,Q)\geq\alpha(P,R)+\alpha(R,Q)=i( \tilde\lambda_+,\sigma)$.

Consider now the inductive step.
Let $\sigma'$ be the segment of $\sigma$ between $x$ and $x_{n-1}$ and $R$ be the 
plane extending the face containing both $x_{n-1}$ and $y$. 
It follows from the inductive hypothesis that
$\alpha(P, R)\geq i(\sigma',\tilde\lambda_+)$.
Still applying Lemma 6.15 of \cite{benedetti-bonsante} to the planes
$P, R, Q$, we obtain that $\alpha(P,Q)\geq i(\sigma',\tilde\lambda_+)+\alpha(R,Q)$.
Notice that the last term is $i(\sigma,\tilde\lambda_+)$ so (\ref{prop1:eq}) 
is proved.

Let us prove now (\ref{prop2:eq}). 
Again in this case we can suppose that $y$ lies on a bending line.
We use again an induction on the number of points $x_1=x,\ldots, x_n=y$
of $\sigma\cap\tilde\lambda_+$.

\begin{figure} 
\begin{center}
\input convexcore.pstex_t
\end{center}
\caption{The proof of (\ref{prop2:eq}).}\label{fig:2}
\end{figure}

When $n=2$, let $F_R$ be the face containing $x$ and $y$ and let $R$
be the plane extending it.  We consider the following surfaces with
boundary: $F_P=P\cap I^-(Q)\cap I^+(R)$, $F_Q=Q\cap I^-(P)\cap I^+(R)$
(see Figure \ref{fig:2}).  Let $\Delta'$ be the surface obtained by
replacing in $\Delta_+$ the face $F_R$ with $F_P\cup F_Q$.

Let $\chi$ be the geodesic segment connecting $x$ to $y$ in $F_P\cup
F_Q$.  Then $\sigma$ and $\chi$ determine a time-like plane in $AdS_3$,
isometric to $AdS_2$.  In $AdS_2$, space-like geodesics are maximizing
length among space-like curves with the same endpoints (this can be
checked easily by writing the AdS metric as a warped product,
$-dt^2+\cos^2(t)dx^2$, with the geodesic segment contained in the line
where $t=0$ -- see \cite{beem} for a more general viewpoint). 
So $\chi$ is shorter than $\sigma$.  In other terms, the
geodesic in $\Delta'$ connecting $x$ to $y$ is shorter than $\sigma$, and
(\ref{prop2:eq}) easily follows.

Consider now the case $n>2$.
Let $R$ be the plane extending the face $F_R$ containing $y$ and $x_{n-1}$ and
$P'$ be the plane extending the  face containing $x_{n-2}$ and
$x_{n-1}$ (see Figure \ref{fig:2}). Let us set
\[
F_{P'}=P'\cap I^-(Q)\cap I^+(R)\,,\quad F_Q=Q\cap I^-(P')\cap I^+(R)~.
\]
Consider the bent surface $\Delta'$ obtained by replacing in
$\Delta_+$ the face $F_R$ by $(F_{P'}\cup F_Q)$.  Let $\sigma'$ be the
geodesic path of $\Delta'$ joining $x$ to $y$. We claim that $\sigma'$ is
shorter than $\sigma$.  Indeed consider the path $\chi$ in $\Delta'$ that is
the composition of the geodesic $\chi_1$ joining $x$ to $x_{n-1}$ and the
geodesic $\chi_2$ joining $x_{n-1}$ to $y$.  Then $\chi_1$ is the
segment of $\sigma$ with endpoints $x$ and $x_{n-1}$, whereas, as before,
$\chi_2$ is shorter than the segment $[x_{n-1},y]$.  Thus $\chi$ is shorter
than $\sigma$. Since $\sigma'$ is shorter than $\chi$, the claim is proved.

Let $y'$ be the intersection point of $\sigma'$ with $P'\cap Q$ and $\sigma''$
be the segment on $\sigma'$ between $x$ and $y'$.  By the inductive
hypothesis we have that $d_P(x, r)+d_Q(y', r)\leq l(\sigma'')$ so we deduce
that $d_P(x,r)+d_Q(y,r)\leq l(\sigma')\leq l(\sigma)$.
\end{proof}

\begin{sublemma}\label{sublemma3}
Let $T$ be a triangle in $AdS_2$ formed by two spacelike rays
$l_1,l_2$ starting from a point $p$ and a complete spacelike geodesic
$l_0$ joinining the ideal end-points of $l_1$ and $l_2$.  Let $q$ be
the point on $l_1$ whose distance from $p$ is $1$ and $\tau$ be the
time-like ray through $q$ orthogonal to $l_1$. If we put $\theta_0=l(\tau\cap T)$,
then
\[
\theta_0=\arctan\left(e\frac{\sh\kappa}{1+\ch\kappa}\right)~,
\]
where $\kappa$ is the angle between $l_1$ and $l_2$.
\end{sublemma}

\begin{figure} 
\begin{center}
\input triangle.pstex_t
\end{center}
\caption{Sublemma \ref{sublemma3}.}\label{fig:triangle}
\end{figure}

\begin{proof}
We identify $AdS_2$ with the quadric in $\mathbb R^3$ given by the
equation $x^2-y^2-z^2=-1$ and equipped with the metric induced by the
form $dx^2-dy^2-dz^2$.  Under such identification the boundary of
$AdS_2$ is identified with the set of light-like vectors up to
multiplication by a positive factor (see Section 2.2).

Under this identification we can suppose that $p=(0,1,0)$ and that the
ideal end-point of $l_1$ is the class of $x_1=(1,1,0)$.  Imposing that
the angle between $l_1$ and $l_2$ is $\kappa$ we deduce that the ideal
end-point of $l_2$ is the class of $x_2=(-\ch\kappa, 1,-\sh\kappa)$
(see Figure \ref{fig:triangle}).
In particular the line $l_0$ is the intersection of $AdS_2$ with the
linear plane generated by $x_1$ and $x_2$. So its equation is
\begin{equation}
(x-y)\sh\kappa- z(1+\ch\kappa)=0\,.
\end{equation}
On the other hand, the coordinates of $q$ are $(\sh 1,\ch 1,0)$, so a
parameterization of $\tau$ is given by $(\cos(\theta) \sh (1),
\cos(\theta) \ch (1),-\sin\theta)$.  Imposing that $\tau(\theta_0)\in
l_0$ we get the result.
\end{proof}

We can now state the estimate we need.

\begin{lemma}\label{lm:curve}
There exist $\kappa_0, \epsilon_1 >0$ as follows.  Let $M$ be an AdS
globally hyperbolic spacetime homeomorphic to $S\times \R$ such that
the upper lamination $\lambda_+$ is a weighted curve $(c,w)$.  Let
$x\in c$ and let $\tau$ be the past-directed geodesic starting from
$x$ and orthogonal to some face $F$.

Given any geodesic $\sigma$ of $\partial_+ C(M)$ that intersects $F$ only
at its end-point $x$ and such that its length is equal $1$ we have
\begin{equation}\label{deep:eq}
   l(\tau\cap C( M))\geq \epsilon_1\min(\kappa_0, i(\sigma,\lambda_+)).
\end{equation}
\end{lemma}

\begin{figure}
\begin{center}
\input bent.pstex_t
\end{center}
\caption{Proof of Lemma \ref{lm:curve}.}
\label{fg:3}
\end{figure}

\begin{proof}[Proof of Lemma \ref{lm:curve}]
We lift the problem to the universal covering.  Let $\Delta_+$ be the
upper boundary of the lift of $C(\tilde M)$ in $\tilde{AdS}_3$, and
let $\xb$ be a lift of $x$.  We consider the lifts $\taub$ and $\bar\sigma$
of $\tau$ and $\sigma$ from $\xb$.  Inequality (\ref{deep:eq}) is
equivalent to $l(\taub\cap C(\tilde M))\geq \epsilon_0\min(\kappa_0,
i(\bar\sigma,\tilde\lambda_+))$

Let $y$ be the end-point of $\bar\sigma$ and consider the plane $P$ through
$\xb$ orthogonal to $\taub$ (that is, the support plane of $\Delta_+$
which extends the face that does not meet $\bar\sigma$) and the plane $Q$
extending the face of $\Delta_+$ that contains $y$ (if $y$ lies on a
bending line we choose the face that does not meet $\bar\sigma$).  By
Sublemma \ref{sublemma}, the distance on $P$ between $\xb$ and the
line $r=P\cap Q$ is less than $1$. So there exists a geodesic
$r'\subset P\cap I^+(Q)$ at distance $1$ from $x$.  Let us consider
the plane $Q'$ containing $r'$ and such that the angle between $P$ and
$Q'$ is $\kappa=i(\bar\sigma,\tilde\lambda_+)$ (see Figure \ref{fg:3}).

By Sublemma \ref{sublemma}, the angle between $P$ and $Q$ is larger
than $\kappa$.  This implies that $Q'$ cannot meet the half-plane
$Q\cap I^-(P)$.  Otherwise the surface obtained as the union of $P\cap
I^-(Q)$, $Q\cap I^-(P)\cap I^-(Q')$ and $Q'\cap I^-(Q)$ would be a 
space-like surface bent along the geodesics $P\cap Q$ and $Q\cap Q'$,
contradicting Equation (\ref{prop1:eq}) in Sublemma \ref{sublemma}.

In particular the surface $\Delta=(P\cap I^-(Q'))\cup (Q'\cap I^-(P))$
is contained in the future of $\Delta_+$. Let $C(\Delta)$ be the convex
hull of $\Delta$.  Clearly the past boundary of $C(\Delta)$ is
contained in the future of the past boundary of $C(\tilde M)$. It
follows that $l(\taub\cap C(\tilde M))$ is larger than $l(\taub\cap
C(\Delta))$.

Now consider the timelike plane $\Pi$ through $\xb$ that is orthogonal
to the bending line $r'$ of $C(\Delta)$.  Notice that $C(\Delta)\cap \Pi$ is
a convex set whose upper boundary is made of two geodesic rays $l_P$
and $l_{Q'}$ meeting at some point $z$. Let $T$ be the triangle in $\Pi$ bounded by
$l_P$, $l_{Q'}$ and by the geodesic joining the ideal end-point of $l_P$
and $l_{Q'}$. $T$ is contained in $C(\Delta)$.  The vertex $z$ is the
point on $r'$ realizing the distance from $\xb$, so the distance
between $\xb$ and $z$ is $1$. Moreover, the angle between $l_P$ and
$l_{Q'}$ is $\kappa$. By Sublemma \ref{sublemma3} the length $l(\taub\cap
T)$ is larger than some constant depending only on $\kappa$ and
proportional to $\kappa$ as $\kappa\rightarrow 0$.  Since $l(\taub\cap
C(\tilde M))\geq l(\taub\cap C(\Delta))\geq l(\taub\cap T)$ the estimate
is proved.
\end{proof}

\begin{remark} \label{rk:pi3}
Choosing $\kappa_0$ sufficiently small, we can suppose that
$\epsilon_1\kappa_0\leq \pi/3$. We make this hypothesis in the rest of
this section.
\end{remark}

\subsection{Pleated surfaces in $AdS_3$.}

We will be using a technical statement from 2-dimensional hyperbolic geometry. The proof
can be found in Appendix \ref{ap:recurrence}.
Consider any  hyperbolic metric $g$ on $S$, and a closed oriented 
geodesic $c$ for $g$. 
Given a point $x\in c$, we consider the compact
geodesic segment $\sigma^l_x$ (resp. $\sigma^r_x$) of length $1$ starting
from $x$ in the direction orthogonal to $c$ towards the left (resp. right). 
We denote by $n^l(x)$ (resp. $n^r(x)$) the number of intersections of $\sigma^l_x$
(resp. $\sigma^r_x$) with $c$, including $x$. Let $\lambda$
be the measured geodesic lamination with support $c$ and weight $w$.

\begin{lemma} \label{lm:inter}
There exists $\beta_0>0$
(depending only on the genus of $S$) such that
$$ l_g(\{ x\in c ~ | ~ i(\sigma^r_x, \lambda)\leq \beta_0 l_g(\lambda)\})
\leq l_g(c)/2~. $$
\end{lemma}

In other terms, if the length of $c$ is much larger than
$\frac{1}{\beta_0}$ then,
in at least half of $c$, the orthogonal segment of length
$1$ on either side of $c$ intersects $c$  
many times. 

Lemma \ref{lm:inter} and Lemma \ref{lm:curve}, taken together, lead to
a lower bound on the area of a totally geodesic time-like strip in the convex core of a
globally hyperbolic AdS manifold. Here we fix a globally hyperbolic
manifold $M$ and we call $C(M)$ its convex core, $\lambda_+$ and
$\lambda_-$ the upper and lower measured bending laminations, and
$m_+, m_-$ the upper and lower induced metrics on the boundary of
$C(M)$. We use the constants $\kappa_0$ and $\epsilon_1$ appearing 
in Lemma \ref{lm:curve}.

\begin{lemma}\label{lm:area}
Suppose that the support of $\lambda_+$ is a closed geodesic curve $c$, on which
we choose an orientation. 
Let $A$ be the time-like totally geodesic annulus in $C(M)$, with boundary 
contained in $\dr C(M)$, such that $c=A\cap \dr_+C(M)$, and 
which is orthogonal to the face on the left of $c$.
Then the area of $A$ is bounded from below:
$$ \area(A) \geq \frac{l_{m_+}(c)}{4} 
\epsilon_1 \min(\kappa_0, \beta_0 l_{m_+}(\lambda_+))~. $$
\end{lemma}

\begin{proof}
Let $w>0$ be the weight of $c$ for $\lambda_+$, and let 
$c_i\subset c$ be the subset of $c$ defined as
$$ c_i = \{ x\in c~ | ~i(\sigma^r_x, \lambda_+)\geq \beta_0 l_{m_+}(\lambda_+)\}~. $$
Lemma \ref{lm:inter} indicates that the length of $c_i$
is at least 
\begin{equation} \label{eq:ci}
l_{m_+}(c_i) \geq l_{m_+}(c)/2~.
\end{equation}

Let $x\in c_i$,  and let $\tau_x$ be the past-directed geodesic 
segment on $A$ starting from $x$ in the direction orthogonal to $c$.
 Because $x\in c_i$  the intersection of $\sigma^r_x$ with the lamination 
 $\lambda_+$ is larger than $\beta_0l_{m_+}(\lambda_+)$.
Applying Lemma \ref{lm:curve}, we deduce that $l(\tau_x)\geq L_1$ 
where
$$ L_1:=\epsilon_1\min(\kappa_0,\beta_0l_{m_+}(\lambda_+))~. $$ 
Note that, by the choice made in Remark \ref{rk:pi3},
\begin{equation} \label{eq:pi3}
   L_1\leq \pi/3~.
\end{equation}

Consider the map
$$ \begin{array}{cccc}
\tau: & c_i\times [0,L_1] & \rightarrow & A \\
& (x,s) & \mapsto & \tau_x(s)~.
\end{array} $$
The Jacobian of $\tau$ at $(x,s)$ is equal to the norm at $\tau_x(s)$
of the Jacobi field along $\tau_x$ which is orthogonal to $\tau_x$
and of unit norm at $s=0$. 
Recall that orthogonal Jacobi fields along time-like geodesics in
$AdS_3$ behave as $\cos(s)$ (it follows from the fact that the 
curvature is $-1$, see e.g. \cite[4.4]{HE}). 
So the Jacobian of $\tau$ at $(x,s)$ is equal to $\cos(s)$,  
and therefore at least $1/2$ by Equation (\ref{eq:pi3}).

As a consequence, 
$$ Area(A)\geq Area(\tau(c_i\times [0,L_1]))\geq \frac{l_{m_+}(c_i)L_1}{2}~.$$
The result therefore follows from
the definition of $L_1$ and from Equation (\ref{eq:ci}).
\end{proof}

\subsection{Proof of the main estimate.}

The proof of Proposition \ref{pr:main} follows from Lemma \ref{lm:area} and
from the following basic statement on Lorentz geometry. 

\begin{lemma}
Let $\Pi$ be a time-like plane in $AdS_3$, and let $P,Q$ be two
space-like planes, such that $\Pi$, $P$ and $Q$ meet exactly at one point.  
Then the angle in $\Pi$ between $\Pi\cap P$ and $\Pi\cap Q$ is smaller than the
angle between $P$ and $Q$.
\end{lemma}

\begin{proof}
Let $x$ be the intersection point of $\Pi, P$ and $Q$, and let
$H_x\subset T_xAdS_3$ be the surface containing the future-oriented 
unit timelike vectors. Clearly $H_x$ is isometric to the hyperbolic
plane. The unit future-pointing vectors $N_P$ and $N_Q$
that are orthogonal to $P$ and $Q$ respectively lie in $H_x$.
Equation (\ref{eq:angle}) shows that the angle between $P$ and $Q$ is equal to
the hyperbolic distance between $N_P$ and $N_Q$ in $H_x$.

On the other hand the set of future-pointing unit vectors tangent to
$\Pi$ at $x$ is the geodesic $l= H_x\cap T_x\Pi$ of $H_x$.  The
future-pointing unit vectors in $T_x\Pi$ that are orthogonal to
$P\cap\Pi$ and $Q\cap \Pi$ are the orthogonal projections of $N_P$ and
$N_Q$ on $l$.

The result therefore follows from the fact that the orthogonal 
projection on a hyperbolic geodesic is distance-decreasing.
\end{proof}

Let $c_-$ be the lower boundary of $A$. 
Although $c_-$ is not geodesic, it is not difficult to see that the geometric 
intersection of $\lambda_-$ with $c_-$ is equal to $i(c,\lambda_-)$.
The previous lemma, along with an approximation of $c_-$ by a sequence of
polygonal curves, shows that the geometric intersection of $\lambda_-$ with
$c_-$ is at least equal to the geodesic
curvature of $c_-$ as boundary of $A$,  which by the Gauss-Bonnet theorem is equal to
the area of $A$ (see \cite{avez:gb,chern:gb}). 
So we obtain the following estimate.

\begin{cor} \label{cr:intersection}
Under the hypothesis of Lemma \ref{lm:area}, if $\lambda_-$ is the measured
bending lamination on the lower boundary of $C(M)$, then
$i(c,\lambda_-) \geq \area(A)$, so that 
$$ i(\lambda_+,\lambda_-)\geq \frac{l_{m_+}(\lambda_+)}{4} 
\epsilon_1 \min(\kappa_0, \beta_0 l_{m_+}(\lambda_+))~. $$
\end{cor}

We can prove now Proposition \ref{pr:main}.

\begin{proof}[Proof of Proposition \ref{pr:main}]
We first  find constants $\epsilon_0, h_0$ that work assuming that
$\lambda$ is a weighted curve $(c,w)$. By a
density argument we then conclude that those constants work for every
measured geodesic lamination.

According to Theorem \ref{tm:mess}, there is a unique GHMC AdS manifold 
$M$ for which the left representation is $g$, the upper bending lamination
of the convex core is $\lambda_+:=\lambda/2$, and the lower bending lamination of
the convex core is $\lambda_-:=\mu/2$. Then the upper and lower induced metrics
on the boundary of the convex core are
$$ m_+ = E_r^{\lambda/2}(g)~,~~ m_- = E_l^{\mu/2}(g)~. $$
We can now apply Corollary \ref{cr:intersection}, which shows that 
$$ i(\lambda_+,\lambda_-) \geq \frac{l_{m_+}(\lambda_+)}{4} 
\epsilon_1 \min(\kappa_0, \beta_0 l_{m_+}(\lambda_+))~. $$
But $g=E_l^{\lambda_+}(m_+)$, so that $l_g(\lambda_+)=l_{m_+}(\lambda_+)$,
and it follows that 
$$ i(\lambda,\mu) \geq \frac{l_g(\lambda)}{2}
\epsilon_1 \min(\kappa_0, \beta_0 l_g(\lambda)/2)~. $$

So the constants
\[
    \epsilon_0=\beta_0\epsilon_1/4~,~~ h_0=2\kappa_0/\beta_0
\]
work.

Consider now the general case.  We can consider a sequence of weighted
curves $\lambda_n$ converging to $\lambda$ such that $\lambda_n=(c_n,
w_n)$.

Let $\mu_n$ be the lamination such that $E^{\lambda_n}_l(g)=E^{\mu_n}_r(g)$.
Notice that $(\mu_n)$ converges to $\mu$ as $n\rightarrow+\infty$. Indeed,
$\mu_n$ is the right factor of $\Phi'(g, E^{t\lambda_n}(g))$.

Since the inequality
\[
   i(\lambda_n,\mu_n)>\epsilon_0l_g(\lambda_n)\min(h_0, l_g(\lambda_n))
\]
holds for every $n$ and $i(\cdot,\cdot)$ is a continuous function of $\cML^2_S$,
passing to the limit we get the estimate for $\lambda$ and $\mu$.
\end{proof}

\section{Compactness} \label{sc:5}

\subsection{Limits of fixed points of compositions of earthquakes}

The rather technical result in the previous section can be used to prove
the following compactness statement, which is a key point in the proof of
the main result. 

\begin{prop} \label{pr:proper}
Let $(\lambda_n)_{n\in \N},(\mu_n)_{n\in \N}$ be two sequences of
measured laminations on $S$, converging respectively to $\lambda,\mu$. 
Suppose that $\lambda$ and $\mu$ fill $S$. Let $(g_n)_{n\in \N}$ be 
a sequence of hyperbolic metrics on $S$ such that, for all $n\in \N$,
$E_l^{\lambda_n}(g_n)=E_r^{\mu_n}(g_n)$. Then, after extracting a subsequence,
$(g_n)$ converges to a limit $g\in \cT_S$, and $E_l^\lambda(g)=E_r^\mu(g)$.
\end{prop}

We recall a well-known result needed in the proof.

\begin{lemma}[Kerckhoff \cite{kerckhoff}] \label{lm:fill}
Let $\lambda,\mu\in \cML_S$ be two measured laminations that fill a surface. 
The function 
$$
\begin{array}{cccc}
L_{\lambda,\mu}: & \cT_S & \rightarrow & \R \\
& g & \mapsto & l_g(\lambda)+l_g(\mu) 
\end{array}
$$
is proper and convex along earthquakes paths.
\end{lemma}

\begin{proof}[Proof of Proposition \ref{pr:proper}]
Since $i(\cdot, \cdot)$ is continuous, there exists
a constant $C_1>0$ such that for all $n\in \N$, 
$i(\lambda_n,\mu_n)\leq C_1$. So it follows from Corollary \ref{cr:obvious} 
that both $l_{g_n}(\lambda_n)$ and $l_{g_n}(\mu_n)$
are bounded by a constant $C_2$. It follows from 
Lemma \ref{lm:fill} that, for all $n\in \N$, 
$$ g_n\in K_n := L^{-1}_{\lambda_n,\mu_n}([0,2C_2])~. $$

Since by Lemma \ref{lm:fill} the subsets $K_n$ are compact and convex for earthquake paths,
$K_n\rightarrow K:=L^{-1}_{\lambda,\mu}([0,2C_2])$.
So $g_n$ remains in a compact subset of $\cT_S$ and, after
taking a subsequence, $(g_n)$ converges to a limit $g$. 

The fact that $E_l^\lambda(g)=E_r^{\mu}(g)$ is clear since $E_l$ and
$E_r$ are continuous functions of both arguments.
\end{proof}

%% file: convexcore.pstex_t
\begin{picture}(0,0)%
\epsfig{file=convexcore.pstex}%
\end{picture}%
\setlength{\unitlength}{1224sp}%
\begingroup\makeatletter\ifx\SetFigFont\undefined%
\gdef\SetFigFont#1#2#3#4#5{%
  \reset@font\fontsize{#1}{#2pt}%
  \fontfamily{#3}\fontseries{#4}\fontshape{#5}%
  \selectfont}%
\fi\endgroup%
\begin{picture}(19742,5314)(632,-5374)
\put(15496,-796){\makebox(0,0)[lb]{\smash{{\SetFigFont{5}{6.0}{\rmdefault}{\mddefault}{\updefault}{\color[rgb]{0,0,0}$r$}%
}}}}
\put(16496,-1796){\makebox(0,0)[lb]{\smash{{\SetFigFont{5}{6.0}{\rmdefault}{\mddefault}{\updefault}{\color[rgb]{0,0,0}$Q$}%
}}}}
\put(12500,-2100){\makebox(0,0)[lb]{\smash{{\SetFigFont{5}{6.0}{\rmdefault}{\mddefault}{\updefault}{\color[rgb]{0,0,0}$P$}%
}}}}
\put(18151,-2656){\makebox(0,0)[lb]{\smash{{\SetFigFont{5}{6.0}{\rmdefault}{\mddefault}{\updefault}{\color[rgb]{0,0,0}$F_Q$}%
}}}}
\put(18951,-2756){\makebox(0,0)[lb]{\smash{{\SetFigFont{7}{8.0}{\rmdefault}{\mddefault}{\updefault}{\color[rgb]{0,0,0}$y$}%
}}}}
\put(17206,-3166){\makebox(0,0)[lb]{\smash{{\SetFigFont{5}{6.0}{\rmdefault}{\mddefault}{\updefault}{\color[rgb]{0,0,0}$F_R$}%
}}}}
\put(16696,-2701){\makebox(0,0)[lb]{\smash{{\SetFigFont{5}{6.0}{\rmdefault}{\mddefault}{\updefault}{\color[rgb]{0,0,0}$F_{P'}$}%
}}}}
\put(15696,-2501){\makebox(0,0)[lb]{\smash{{\SetFigFont{5}{6.0}{\rmdefault}{\mddefault}{\updefault}{\color[rgb]{0,0,0}$x_{n-1}$}%
}}}}
\put(3766,-3631){\makebox(0,0)[lb]{\smash{{\SetFigFont{5}{6.0}{\rmdefault}{\mddefault}{\updefault}{\color[rgb]{0,0,0}$F_R$}%
}}}}
\put(5026,-2746){\makebox(0,0)[lb]{\smash{{\SetFigFont{5}{6.0}{\rmdefault}{\mddefault}{\updefault}{\color[rgb]{0,0,0}$F_Q$}%
}}}}
\put(4201,-1846){\makebox(0,0)[lb]{\smash{{\SetFigFont{5}{6.0}{\rmdefault}{\mddefault}{\updefault}{\color[rgb]{0,0,0}$r$}%
}}}}
\put(3106,-2551){\makebox(0,0)[lb]{\smash{{\SetFigFont{5}{6.0}{\rmdefault}{\mddefault}{\updefault}{\color[rgb]{0,0,0}$F_P$}%
}}}}
\end{picture}%

%% file: triangle.pstex_t
\begin{picture}(0,0)%
\epsfig{file=triangle.pstex}%
\end{picture}%
\setlength{\unitlength}{3158sp}%
\begingroup\makeatletter\ifx\SetFigFont\undefined%
\gdef\SetFigFont#1#2#3#4#5{%
  \reset@font\fontsize{#1}{#2pt}%
  \fontfamily{#3}\fontseries{#4}\fontshape{#5}%
  \selectfont}%
\fi\endgroup%
\begin{picture}(8101,3699)(1051,-5173)
\put(8476,-2761){\makebox(0,0)[lb]{\smash{{\SetFigFont{10}{12.0}{\rmdefault}{\mddefault}{\updefault}{\color[rgb]{0,0,0}$[1,1,0]$}%
}}}}
\put(5776,-2611){\makebox(0,0)[lb]{\smash{{\SetFigFont{10}{12.0}{\rmdefault}{\mddefault}{\updefault}{\color[rgb]{0,0,0}$(1,0,0)$}%
}}}}
\put(3826,-2986){\makebox(0,0)[lb]{\smash{{\SetFigFont{10}{12.0}{\rmdefault}{\mddefault}{\updefault}{\color[rgb]{0,0,0}$(-\ch\kappa,0,-\sh\kappa)$}%
}}}}
\put(5326,-2611){\makebox(0,0)[lb]{\smash{{\SetFigFont{10}{12.0}{\rmdefault}{\mddefault}{\updefault}{\color[rgb]{0,0,0}$p$}%
}}}}
\put(6601,-2611){\makebox(0,0)[lb]{\smash{{\SetFigFont{10}{12.0}{\rmdefault}{\mddefault}{\updefault}{\color[rgb]{0,0,0}$q$}%
}}}}
\put(6826,-3061){\makebox(0,0)[lb]{\smash{{\SetFigFont{10}{12.0}{\rmdefault}{\mddefault}{\updefault}{\color[rgb]{0,0,0}$\tau$}%
}}}}
\put(5176,-3811){\makebox(0,0)[lb]{\smash{{\SetFigFont{10}{12.0}{\rmdefault}{\mddefault}{\updefault}{\color[rgb]{0,0,0}$l_0$}%
}}}}
\put(1051,-4561){\makebox(0,0)[lb]{\smash{{\SetFigFont{10}{12.0}{\rmdefault}{\mddefault}{\updefault}{\color[rgb]{0,0,0}$[-\ch\kappa,1,-\sh\kappa]$}%
}}}}
\end{picture}%

%% file: bent.pstex_t
\begin{picture}(0,0)%
\epsfig{file=bent.pstex}%
\end{picture}%
\setlength{\unitlength}{2289sp}%
\begingroup\makeatletter\ifx\SetFigFont\undefined%
\gdef\SetFigFont#1#2#3#4#5{%
  \reset@font\fontsize{#1}{#2pt}%
  \fontfamily{#3}\fontseries{#4}\fontshape{#5}%
  \selectfont}%
\fi\endgroup%
\begin{picture}(10440,5664)(2098,-4573)
\put(2896,-3451){\makebox(0,0)[lb]{\smash{{\SetFigFont{7}{8.4}{\rmdefault}{\mddefault}{\updefault}{\color[rgb]{0,0,0}$P$}%
}}}}
\put(7996,509){\makebox(0,0)[lb]{\smash{{\SetFigFont{7}{8.4}{\rmdefault}{\mddefault}{\updefault}{\color[rgb]{0,0,0}$z$}%
}}}}
\put(9091,-2131){\makebox(0,0)[lb]{\smash{{\SetFigFont{7}{8.4}{\rmdefault}{\mddefault}{\updefault}{\color[rgb]{0,0,0}$Q$}%
}}}}
\put(11386,-136){\makebox(0,0)[lb]{\smash{{\SetFigFont{7}{8.4}{\rmdefault}{\mddefault}{\updefault}{\color[rgb]{0,0,0}$Q'$}%
}}}}
\put(8896, 59){\makebox(0,0)[lb]{\smash{{\SetFigFont{7}{8.4}{\rmdefault}{\mddefault}{\updefault}{\color[rgb]{0,0,0}$r'$}%
}}}}
\put(6796,-1111){\makebox(0,0)[lb]{\smash{{\SetFigFont{7}{8.4}{\rmdefault}{\mddefault}{\updefault}{\color[rgb]{0,0,0}$r$}%
}}}}
\put(3976,-2611){\makebox(0,0)[lb]{\smash{{\SetFigFont{7}{8.4}{\rmdefault}{\mddefault}{\updefault}{\color[rgb]{0,0,0}$\bar x$}%
}}}}
\end{picture}%

%% file: fb3.tex
\section{Proofs of the main results} \label{sc:6}

In this section we combine results of Sections \ref{sc:3} and \ref{sc:5} 
to prove Theorems \ref{tm:1}, \ref{tm:2} and \ref{tm:4}. 

Consider the map $\Phi'$ described at the beginning of Section
\ref{sc:2}.  Proposition~\ref{pr:proper} precisely states that this
map is proper.  In particular the degree of $\Phi'$ can be defined.
The main idea of the argument is to show that the degree of $\Phi'$ is
$1$. As a consequence, we will deduce that $\Phi'$ is surjective, and
this will prove Theorem \ref{tm:1} and Theorem \ref{tm:2}.

\subsection{Proof of Theorem \ref{tm:1} and Theorem \ref{tm:2}}

Let us consider the set
\[
   X=\{(u,v)\in\mathcal T^2| \Phi'^{-1}\Phi'(u,v)=\{(u,v)\}\}~.
\]

\begin{prop}\label{open:prop}
$\hat X=\Delta\cup X$ is a neighbourhood of $\Delta$ in $\mathcal T^2$.
\end{prop}

\begin{proof}
By contradiction, suppose there exists a sequence $(u_k,v_k)\notin X$
converging to $(u,u)\in\Delta$.  By Proposition \ref{pr:U}, 
there exists another sequence
$(u'_k, v'_k)\in\mathcal T$ such that
\[
   \Phi'(u'_k, v'_k)=\Phi'(u_k, v_k)~.
\]

By Lemma \ref{tec2:lem} there is an infinitesimal sequence of positive
numbers $t_k$ such that
\[
   \Phi'(u'_k,v'_k)=\Phi'(u_k, v_k)=(t_k\lambda_k, t_k\mu_k)
\]
with $\{(\lambda_k,\mu_k)\}$ running in some compact set of
$\mathcal F\mathcal M\mathcal L_S$. Taking a subsequence we
can suppose that $(\lambda_k,\mu_k)\rightarrow(\lambda,\mu)$ and
 $u$ is the point  $k_0(\lambda,\mu)$ where $e^\lambda_l$ and $e^\mu_r$ meet.

By Proposition \ref{pr:main} we have
\[
    t_k^2 i(\lambda_k,\mu_k)\geq \epsilon_0 t_k^2 l_{u'_k}(\lambda_k)^2
\]
so we get that $l_{u'_k}(\lambda_k)$ is bounded by some constant
independent of $k$.  Analogously we deduce that $l_{u'_k}(\mu_k)$ is
bounded.  Since $(\lambda_k,\mu_k)$ runs in a compact set of 
$\mathcal F\mathcal M\mathcal L_S$, the functions
\[
     l(\lambda_k)+l(\mu_k)
\]
are uniformly proper  by Lemma \ref{lm:fill}.  

So we deduce that $(u'_k)$ runs in some compact set of $\mathcal T$ (and
analogously for $(v'_k)$). Moreover, by Lemmas \ref{tec1:lem} and
\ref{tec2:lem}, any convergent subsequence of $(u'_k)$ (resp. $(v'_k)$)
must converge to $k_0(\lambda,\mu)=u$. Thus we deduce that
$(u'_k,v'_k)\rightarrow (u,u)$.

By Corollary \ref{cr:loc}, there is a neighbourhood $U$ of
$(u,u)\in\mathcal T^2$ such that the map $\Phi'|_{U\setminus\Delta}$ is
a homeomorphism onto an open set in $\mathcal F\mathcal M\mathcal L_S$.
But for $k>>0$ both $(u_k,v_k)$ and $(u'_k, v'_k)$ lie in
$U$ and this is a contradiction.
\end{proof}

\begin{cor} \label{cr:open}
There is an open set $V\in\mathcal F\mathcal M\mathcal L_S$ such that
the restriction of $\Phi'$ to $\Phi'^{-1}(V)$ is a homeomorphism onto
$V$.
\end{cor}

\begin{proof}
There is an open neighbourhood $U$ of $\Delta$ contained in $X$.  Let
us put $V=\Phi'(U)$.  Clearly the restriction of $\Phi'$ on $U$ is
injective so $V$ is an open set.  By definition of $X$ 
the inverse image of any point $x\in V$ consists only of one point.  Thus
$\Phi'^{-1}(V)=U$.
\end{proof}

\begin{cor}
The degree of the map $\Phi'$ is $1$.  In particular, $\Phi'$ is surjective.
\end{cor}

\begin{proof}
This follows from Corollary \ref{cr:open}, since a continuous (proper)
map which restricts to a homeomorphism from an open subset to its image
has degree one, see \cite{hirsch}.
\end{proof}

Since $\Phi'$ has degree one, it is surjective, and this proves Theorem
\ref{tm:1}.  We have already mentioned that Theorem \ref{tm:2} is
equivalent to Theorem \ref{tm:1}.


\subsection{Proof of Theorem \ref{tm:4}}

Let  $(\lambda,\mu)\in\mathcal F\mathcal M\mathcal L_S$. 
We have to prove that for $t$ sufficiently small $(t\lambda, t\mu)$
are uniquely realized as the bending laminations of the convex core
of a GHMC AdS spacetime. By Equation (\ref{mm:eq}), this is equivalent to showing
that $\Phi'^{-1}(t\lambda, t\mu)$ contains exactly one point for $t$ small.

We consider the right inverse of $\Phi'$, $\sigma:V\rightarrow \mathcal
T\times\mathcal T$, defined in Corollary \ref{cr:loc}. There is
$\epsilon>0$ such that $(t\lambda, t\mu)\in V$ for $t<\epsilon$.  In
particular $\sigma(t\lambda, t\mu)$ is in $\Phi'^{-1}(t\lambda,t\mu)$. On
the other hand, $\sigma(t\lambda, t\mu)$ approaches the diagonal
$\Delta$ as $t\rightarrow 0$. By Proposition \ref{open:prop}, there is
$\epsilon'\leq\epsilon$ such that $\Phi'^{-1}(t\lambda,
t\mu)=\{\sigma(t\lambda, t\mu)\}$.

\subsection{Surfaces with cone singularities}

As pointed out in the introduction, the arguments given for the proof
of Theorem \ref{tm:1} and of Theorem \ref{tm:2} can be extended
basically as they are to surfaces with cone singularities of angle
$\theta_i\in (0,\pi)$, and to GHMC AdS manifolds with ``particles'' of
the same angles. Thurston's Earthquake Theorem is then replaced by its
version with ``particles'' as described in \cite{cone}, where the
geometry of GHMC AdS manifolds with particles was also studied.
This leads directly to the proof of Theorem \ref{tm:5} or,
equivalently, Theorem \ref{tm:6}.

%% file: recurrence.tex
\section{Proof of Lemma \ref{lm:inter}}
\label{ap:recurrence}

Recall that we consider a hyperbolic metric $g$ on $S$, and a closed oriented 
geodesic $c$ for $g$. 
Given $x\in c$, $\sigma^l_x$ (resp. $\sigma^r_x$) is the compact geodesic segment of length $1$ starting
from $x$ in the direction orthogonal to $c$ towards the left (resp. right), 
and $n^l(x)$ (resp. $n^r(x)$) is the number of intersections of $\sigma^l_x$
(resp. $\sigma^r_x$) with $c$, including $x$. 
Lemma \ref{lm:inter} is a direct consequence of the following statement.

\begin{lemma} \label{lm:A1}
There exists $\beta_0>0$
(depending only on the genus of $S$) such that
$$ l_g(\{ x\in c~ | ~ n^r(x)\leq \beta_0 l_g(c) \}) \leq l_g(c)/2 $$
for every simple closed geodesic $c$.
\end{lemma}

The proof uses an elementary statement from plane hyperbolic geometry.

\begin{sublemma} \label{slm}
There exists $\gamma_0>0$ as follows. Let $D_0,D_1$ be two disjoint
lines in $\mathbb H^2$, and let $x\in \mathbb H^2$ be in the connected
component of $\mathbb H^2\setminus (D_0\cup D_1)$ having both $D_0$
and $D_1$ in its boundary.  Suppose that $d(x,D_0)\leq \gamma_0$,
$d(x,D_1)\leq \gamma_0$. Then the geodesic segment of length $1$
starting orthogonally from $D_0$ and containing $x$ intersects $D_1$.
\end{sublemma}

\begin{figure}
\begin{center}
\scalebox{0.7}{
\input sublemma2.pstex_t
}
\end{center}
\caption{Sublemma \ref{slm}}.
\end{figure}

\begin{proof}
It is sufficient to find $\gamma_0$ assuming that $D_0$ is a fixed
geodesic and that the end-point on $D_0$ of the minimizing geodesic
segment going from $x$ to $D_0$ is a fixed point
$y_0$. Finally we can assume that $x$ is contained in a fixed half-plane 
$P_0$ bounded by $D_0$.

Let $x_0\in P_0$ be the end-point of the geodesic segment of length
$1$ starting from $y_0$ and orthogonal to $D_0$, and let $l_1, l_2$ be
the two complete geodesics that share an ideal end-point with $D_0$
and pass through $x_0$.

Notice that either $D_1$ meets the segment $[y_0,x_0]$ or its distance
from $y_0$ is larger than the distance $\delta$ of the lines $l_i$
from $y_0$.
Now, taking $\gamma_0=\delta/2$, the distance between $D_1$ and $y_0$
is less than $d(D_1,x)+d(x,y_0)\leq 2\gamma_0=\delta$, so $D_1$
intersects $[y_0,x_0]$.
\end{proof}

\begin{proof}[Proof of Lemma \ref{lm:A1}]
We take here some $\beta_0>0$, and will later see how it has to
be chosen so as to obtain the desired result. Let 
$$ c_i := \{ x\in c~ | ~ n^r(x)\leq \beta_0l_g(c) \}~. $$ Fix
$\gamma_0$ as in Sublemma \ref{slm} and consider the normal
exponential map:
$$
\begin{array}{cccc}
\exp: & c_i\times [0,\gamma_0] & \rightarrow & S \\
& (s,r) & \mapsto & \sigma^r_s(r)~. 
\end{array}
$$
This map is dilating, so it increases the area. 

Moreover, Sublemma \ref{slm} shows that each point $x\in S$ has at
most $n_0$ inverse images in $c_i\times [0,\gamma_0]$, where $n_0$ is
the integer part of $\beta_0 l_g(c)$.  Indeed, suppose that $x$ is the
image of $(y_1,r_1),\ldots, (y_n,r_n)$ .  Let
$\xb,\yb_1,\yb_2,\ldots,\yb_n$ be lifts of $x,y_1,\ldots,y_n$ to the
universal cover $\mathbb H^2$ of $(S,g)$ chosen so that some lift of
$\exp(y_j,[0,\gamma_0])$ contains both $\xb$ and $\yb_j$. Finally let
$D_i$ be the lift of $c$ passing through $y_i$.

For $i\neq j$, $D_i$ and $D_j$ are disjoint. Indeed, since the segment
$[\xb,\yb_i]$ is orthogonal to $D_i$, the lines $D_i$ and $D_j$ cannot
coincide.

Up to changing the indices, we can suppose that there are half-planes
$P_1$ and $P_2$ bounded by $D_1$ and $D_2$ respectively that do not
meet any other leaf $D_i$.  Up to exchanging $D_1$ and $D_2$ we can in
particular suppose that $\xb\notin P_1$.  For $i\geq 2$ either $D_i$
disconnects $D_1$ from $\xb$ or $\xb$ is contained in the region
bounded by $D_1$ and $D_i$.  In the latter case Sublemma \ref{slm} can
be applied since the distance of $\xb$ from $D_1$ and $D_i$ is less
than $\gamma_0$.  In both cases, the segment of length $1$ starting
from $\yb_1$ and passing through $\xb$ --- which is a lift of
$\sigma^r_{y_1}$ --- meets $D_i$.  Since $\sigma^r_{y_1}$ meets $c$ at most
$n_0$ times (including $y_1$), we conclude that $n\leq n_0$.

Since the area of $(S,g)$ is  $2\pi|\chi(S)|$, it follows that
$$ \gamma_0 l_g(c_i) \leq 2\pi n_0 |\chi(S)| 
\leq 2\pi \beta_0 l_g(c)|\chi(S)|~. $$
The result clearly follows, with a choice of
$\beta_0=\gamma_0/(4\pi|\chi(S)|)$.
\end{proof}


%% file: sublemma2.pstex_t
\begin{picture}(0,0)%
\epsfig{file=sublemma2.pstex}%
\end{picture}%
\setlength{\unitlength}{3947sp}%
\begingroup\makeatletter\ifx\SetFigFont\undefined%
\gdef\SetFigFont#1#2#3#4#5{%
  \reset@font\fontsize{#1}{#2pt}%
  \fontfamily{#3}\fontseries{#4}\fontshape{#5}%
  \selectfont}%
\fi\endgroup%
\begin{picture}(6296,6366)(231,-6394)
\put(5296,-4141){\makebox(0,0)[lb]{\smash{{\SetFigFont{12}{14.4}{\rmdefault}{\mddefault}{\updefault}{\color[rgb]{0,0,0}$l_1$}%
}}}}
\put(4486,-3226){\makebox(0,0)[lb]{\smash{{\SetFigFont{12}{14.4}{\rmdefault}{\mddefault}{\updefault}{\color[rgb]{0,0,0}$x_0$}%
}}}}
\put(4921,-2026){\makebox(0,0)[lb]{\smash{{\SetFigFont{12}{14.4}{\rmdefault}{\mddefault}{\updefault}{\color[rgb]{0,0,0}$l_2$}%
}}}}
\put(4081,-1021){\makebox(0,0)[lb]{\smash{{\SetFigFont{12}{14.4}{\rmdefault}{\mddefault}{\updefault}{\color[rgb]{0,0,0}$D_1$}%
}}}}
\put(2761,-1096){\makebox(0,0)[lb]{\smash{{\SetFigFont{12}{14.4}{\rmdefault}{\mddefault}{\updefault}{\color[rgb]{0,0,0}$D_0$}%
}}}}
\put(2896,-3151){\makebox(0,0)[lb]{\smash{{\SetFigFont{12}{14.4}{\rmdefault}{\mddefault}{\updefault}{\color[rgb]{0,0,0}$y_0$}%
}}}}
\put(3601,-3391){\makebox(0,0)[lb]{\smash{{\SetFigFont{12}{14.4}{\rmdefault}{\mddefault}{\updefault}{\color[rgb]{0,0,0}$x$}%
}}}}
\end{picture}%

%% file: fb2.tex
\section{Flat case}

Let $\R^{2,1}$ be the standard $3$-dimensional Minkowski space,
that is $\mathbb R^3$ equipped with the flat Lorentzian metric
$dx_1^2+dx_2^2-dx_3^2$. The isometry group of $\mathbb R^{2,1}$ is the
affine group of transformations whose linear part preserves the
Minkowski product.  Thus we have
\[
   Isom(\R^{2,1})= \mathbb R^3\rtimes O(2,1)
\]
where $O(2,1)$ acts on $\mathbb R^3$ by multiplication.

In this section we will consider the hyperboloid model of
$\mathbb H^2$, that is the set of future-pointing unit time-like vectors
in $\R^{2,1}$. Using this model, we identify the orientation-preserving 
isometry group of $\mathbb H^2$ with the connected component  to the identity 
of $O(2,1)$, which we denote by $SO^+(2,1)$.

Finally we will identify the Teichm\"uller space of $S$ with the set of
conjugacy classes of faithful and discrete representations
$h:\pi_1(S)\rightarrow SO^+(2,1)$.

For any Fuchsian representation
\[
    h:\pi_1(S)\rightarrow SO^+(2,1)
\]
the cones $I^+(0)$ and $I^-(0)$ 
are invariant subsets of $\R^{2,1}$ 
and the quotients
$I^+(0)/h$ and $I^-(0)/h$ are  globally hyperbolic flat spacetimes
homeomorphic to $S\times\mathbb R$.

Now consider an affine deformation of $h$, that is, a representation
$ \rho:\pi_1(S)\rightarrow\mathbb R^3\rtimes SO^+(2,1)$
which is  of the form
\[
   \rho(\alpha)=h(\alpha)+\tau(\alpha)
\]
where $\tau(\alpha)\in\mathbb R^3$ is the translation part.  In
\cite{mess}, Mess showed that there are two maximal convex domains in
$\R^{2,1}$ --- say $\Omega^+(\rho)$ and $\Omega^-(\rho)$ --- such that:
\begin{itemize}
\item They are invariant for the action of $\rho(\pi_1(S))$. The
  restriction of the action of $\rho(\pi_1(S))$ on them is free and
  properly discontinuous.
\item $M^+(\rho)=\Omega^+(\rho)/\rho$ and $M^-(\rho)=\Omega^-(\rho)/\rho$ are
globally hyperbolic spacetimes homeomorphic to $S\times\mathbb R$.
\item $\Omega^+(\rho)$ is complete in the future: if a point $x$ lies
  in $\Omega^+(\rho)$ then every future-directed timelike  path starting at
  $x$ is contained in $\Omega^+(\rho)$.  Analogously $\Omega^-(\rho)$
  is complete in the past.
\end{itemize}

Mess proved that all globally hyperbolic flat spacetimes homeomorphic to
$S\times\mathbb R$ are contained in a spacetime of this form.

\begin{prop}\cite{mess}\label{flat1:prop}
Let $M$ be a time-oriented globally hyperbolic flat spacetime homeomorphic
to $S\times\mathbb R$.  There is an affine deformation of some Fuchsian
representation
\[
  \rho:\pi_1(S)\rightarrow \mathbb R^3\rtimes SO^+(2,1)
\]
such that $M$ isometrically embeds either in $M^+(\rho)$ or in $M^-(\rho)$. 
\end{prop}

We denote by $\mathcal H$ the set of
conjugacy classes of  representations
$\rho:\pi_1(S)\rightarrow\mathbb R^3\rtimes SO^+(2,1)$ whose linear
part is faithful and discrete.  The Fuchsian locus of $\mathcal H$ ---
that we will denote by $\mathcal H_0$ --- corresponds to
representations that are conjugate to Fuchsian representations.

We have a projection map 
$\pi_{\mathcal H}:\mathcal H\rightarrow\mathcal T$,
sending $\rho$ to its linear part.
The fiber over $h\in\mathcal T$ consists of all the elements of the form
\[
    \rho=h+\tau
\]
where $\tau:\pi_1(S)\rightarrow\mathbb R^3$ represents the translation
part.  Then $\rho$ is a representation if and only if 
$\tau$ satisfies the cocycle condition
\[
\tau(\alpha\beta)=\tau(\alpha)+ h(\alpha)\tau(\beta)\,.
\]
Thus $\tau$ is an element of $Z^1_{h}(\pi_1(S),\mathbb
R^3)$, which is the group of $1$-cocycles of $\pi_1(S)$ with values in
$\mathbb R^3$, where the action of $\pi_1(S)$ on $\mathbb R^3$ is
induced by $h$.

Two representations obtained by two cocycles $\tau,\tau'$ are
conjugate if and only if $\tau-\tau'$ is a coboundary.
So the fiber of $h$ is identified to 
the cohomology group $H^1_h(\pi_1(S), \mathbb R^3)$.

It turns out that the fibers of $\pi_{\mathcal H}$ have a natural
structure of vector spaces.  In fact, the map $\pi_{\mathcal
  H}:\mathcal H\rightarrow\mathcal T$ is a vector bundle of rank
$6g-6$ (see \cite{mess}).  Notice that $\mathcal H_0$ is
the image of the zero section.

\subsection{Laminations associated to an affine deformation}

Given some $\rho=h+\tau\in\mathcal H$,
we consider the distance of points in $\Omega^+(\rho)$ from the boundary,
that is, the function defined by the formula 
\[
     \tilde t(x)=\sup_{y\in\partial\Omega^+\cap
       I^-(x)}\big(-\langle(x-y),(x-y)\rangle\big)^{1/2}\,.
\]
This function is $\mathrm C^{1,1}$ and its level sets are achronal
\cite{mess,bonsante}.  It induces a function $t$ on $M^+(\rho)$ and we
consider the level surface $S^+=t^{-1}(1)$.

Notice that when $\rho$ is a Fuchsian representation, then we simply
have $S^+=\mathbb H^2/\rho$.  In the general case, Mess showed that
$S^+$ is obtained by \emph{grafting} the hyperbolic surface
corresponding to the linear part of $\rho$ along a measured geodesic
lamination $\lambda_+$.  
More precisely, the induced metric on $S_+$ is isometric to the 
``grafted metric'' which, if $\lambda_+$ is a weighted multicurve,
is obtained by replacing in the hyperbolic metric associated to
$\rho$ each leave of $\lambda_+$ by a flat cylinder of width
equal to the weight.
Through this section $\lambda_+$ is called the
\emph{upper lamination} of $\rho$.

In the same way, the surface $S^-\subset M^-(\rho)$ of points at
distance $1$ from the boundary is obtained by grafting the hyperbolic
surface corresponding to the linear part of $\rho$ along a measured
geodesic lamination $\lambda_-$, which is called the \emph{lower lamination}
of $\rho$.

We consider the map
\[
    \Phi_0:\mathcal H\rightarrow\mathcal{ML}_S^2
\]
defined by $\Phi_0(\rho)=(\lambda_+,\lambda_-)$.
Notice that if $\rho$ is a Fuchsian representation  
then $\Phi_0(\rho)=(0,0)$.
In this section we study the map $\Phi_0$ outside the Fuchsian locus.
In particular we prove the following theorem.
\begin{theorem}\label{frat:tm}
The map
\[
\Phi_0:\mathcal H\setminus\mathcal H_0\rightarrow\mathcal{ML}_S^2
\]
is a homeomorphism onto $\cFML_S$.
\end{theorem}


In \cite{mess}, it is proved that given a measured geodesic $\lambda$ and
a Fuchsian representation $h$, there is a unique cocycle $\tau_+\in
H^1_h(\pi_1(S),\mathbb R^3)$ such that $\lambda$ is the upper
lamination of $h+\tau_+$.  Analogously there is a unique cocycle
$\tau_-\in H^1_h(\pi(S),\mathbb R^3)$ such that the lower lamination
of $h+\tau_-$ is $\lambda$.

In this way, for a fixed $\lambda$, we determine two sections
$\tau^\lambda_+,\tau^\lambda_-:\mathcal T\rightarrow\mathcal H$
by requiring that $\lambda$ is the upper (resp. lower) lamination of
$h+\tau^\lambda_+(h)$ (resp. $h+\tau^{\lambda}_-(h)$).

These sections are explicitly described in \cite{mess}.  For the
reader's convenience we describe the simple case where the lamination
$\lambda$ is a weighted multicurve, referring to \cite{mess,
  benedetti-bonsante} for details.  Given a Fuchsian representation
$h:\pi_1(S)\rightarrow SO^+(2,1)$, we realize $\lambda$ as a weighted
geodesic multicurve in $\mathbb H^2/h$.  Let $\tilde\lambda$ be 
its lift to $\mathbb H^2$.  Now let us fix a point $x_0\in\mathbb
H^2\setminus\tilde\lambda$.  For any $\alpha\in\pi_1(S)$, we consider
the intersection points --- say $p_1,\ldots, p_N$ --- of the geodesic
segment $[x_0,h(\alpha)(x_0)]$ with $\tilde\lambda$.  Then we define
the following vector in $\mathbb R^3$:
\begin{equation}\label{mess:eq}
    \tau(\alpha)=\sum_{i=1}^N m_iw_i\,,
\end{equation}
where $m_i$ is the weight of the leaf through $p_i$ and $w_i\in\mathbb
R^3$ is the unit vector orthogonal to the leaf through $p_i$ pointing
towards $\alpha(x_0)$. It turns out that $\tau$ is a cocycle and (up
to a coboundary) we have
\begin{equation}\label{mess:eq2}
    \tau^\lambda_+(h)(\alpha)=\tau\qquad\qquad
    \tau^\lambda_-(h)(\alpha)=-\tau\,.
\end{equation}

Notice that if $\Phi_0(h+\tau)=(\lambda_+,\lambda_-)$ then we clearly
have $\tau=\tau^{\lambda_+}_+(h)=\tau^{\lambda_-}_-(h)$. In particular,
the sections $\tau^{\lambda_+}_+$ and $\tau^{\lambda_-}_-$ meet over $h$.
Conversely if for some given $\lambda,\mu\in\mathcal{ML}_S$ there is
$h\in\mathcal T$ such that $\tau^{\lambda}_+(h)=\tau^{\mu}_-(h)$ 
then $\Phi(h+\tau^{\lambda}(h))=(\lambda,\mu)$.  Thus there is
a $1$-to-$1$ correspondence between $\Phi_0^{-1}(\lambda,\mu)$ and
the intersection of $\tau^{\lambda}_+$ and $\tau^{\mu}_-$.

The following proposition combined with Proposition \ref{kerck:prop}
shows that $\tau^{\lambda}_+$ and $\tau^\mu_-$ are disjoint if
$\lambda$ and $\mu$ does not fill, whereas they intersect at exactly one
point otherwise. This shows that the image of $\Phi_0$ is $\cFML_S$
and that $\Phi_0$ is injective.
\begin{prop}
There is a  vector bundle isomorphism
\[
     \xi_*:\mathcal H\rightarrow T\mathcal T_S
\]
such that $\xi_*\circ\tau^\lambda_+=e^\lambda_r$ and
$\xi_*\circ\tau^\lambda_-=e^\lambda_l$.
\end{prop}

\begin{proof}
We consider on $\mathbb R^3$ the Minkowski
vector product: given $x=(x_1,x_2,x_3)$ and $y=(y_1,y_2,y_3)$
it is  defined by
\[
   x\times y=(x_2y_3-x_3y_2,\ x_3y_1-x_1y_3,\  -x_1y_2+x_2y_1)\,.
\]
We refer to \cite{drumm} for the details.
Given $x\in \mathbb R^3$, the linear operator $\xi(x)$ defined by
$\xi(x)(y)=x\times y$ is skew-symmetric with respect to the Minkowski product
$\langle\cdot,\cdot\rangle$, so it lies in the Lie algebra $\mathfrak
o(2,1)$ of $SO^+(2,1)$.  The induced map $\xi:\mathbb R^3\rightarrow
\mathfrak o(2,1)$ is an isomorphism. Moreover we have
\[
   \xi(Ax)=Ad(A)\xi(x)\,.
\]

Given $h\in\mathcal T$ and an element $\tau\in Z^1_h(\pi_1(S),\mathbb
R^3)$ we have that $\tau^*=\xi\circ\tau:\pi_1(S)\rightarrow \mathfrak o(2,1)$
satisfies the cocycle rule
\[
   \tau^*(\alpha\beta)=\tau^*(\alpha)+Ad(h(\alpha))\tau^*(\beta)\,.
\]
In particular, $\tau^*$ represents an infinitesimal deformation of the
representation $h$.  Moreover $\tau^*$ is a trivial deformation if and
only if $\tau$ is a coboundary.  Using the canonical identification
between $T_h\mathcal T$ and $H^1_{Ad\circ h}(\pi_1(S), \mathfrak
o(2,1))$ (see \cite{goldman})
we obtain an isomorphism
\[
   \xi_{*,h}: H^1_h(\pi_1(S),\mathbb R^3)\rightarrow T_h\mathcal T
\]
defined by $\xi_{*,h}([\tau])=[\tau^*]$.

Since $H^1_h(\pi_1(S),\mathbb R^3)$ is the fiber of the
projection $\pi_{\mathcal H}:\mathcal H\rightarrow\mathcal T$ over
$h$, the maps $\xi_{*,h}$ produce a fiber bundle isomorphism
\[
\xi_*:\mathcal H\rightarrow T\mathcal T\,.
\] 

To conclude the proof, we have to prove that
$\xi_*\circ\tau^\lambda_+(h)=e^\lambda_r(h)$
and $\xi_*\circ\tau^\lambda_-(h)=e^\lambda_l(h)$.
Since $\tau^\lambda_-=-\tau^\lambda_+$ and 
$e^\lambda_l=-e^\lambda_r$ it is sufficient to prove only
the first equality.

We will prove that the equality holds when $\lambda$ is a weighted curve.
The general case follows by a simple approximation argument.

We denote by $\tilde\lambda$ the measured geodesic lamination on
$\mathbb H^2$ that projects on the geodesic lamination on $\mathbb
H^2/h$ that realizes $\lambda$.  We fix a point $x_0\in\mathbb
H^2\setminus\tilde\lambda$.  Given $\alpha\in\pi_1(S)$, let
$l_1,\ldots, l_N$ be the leaves of $\tilde\lambda$ meeting the segment
$[x_0,\alpha(x_0)]$ at points $x_1,\ldots, x_N$ respectively.  We
denote by $m_i$ the weight of $l_i$ and by $w_i\in\mathbb R^3$ the
unit spacelike vector orthogonal to $l_i$ and pointing towards
$\alpha(x_0)$.

It can be easily checked that the transformation $\xi(w_i)$ is an
infinitesimal generator of the hyperbolic group of transformations
with axis $l_i$.  More precisely, for $m>0$, a simple computation shows 
that $\exp(m\xi(w_i))$ is the hyperbolic transformation of axis $l_i$ and
translation length equal to $m$ moving $\alpha(x_0)$ to the right of
$x_0$.

Let us put  $h_t=E_r^{t\lambda}(h)$. 
It follows from computations in \cite{epstein-marden} that
\[
    h_t(\alpha)=\exp(tm_1\xi(w_1) )\circ\cdots \circ \exp(tm_N
    \xi(w_N))h(\gamma)\,.
\]
Taking the derivative at $t=0$,
\[
      \frac{ d h_t(\gamma)}{dt}|_0=
      (R_{h(\gamma)})_*\left(\sum
      m_i\xi(w_i)\right)\,.
 \]     

So the cocycle representing the derivative at $0$ of such a deformation
(that is, the cocycle corresponding to $e^\mu_r(h)$ via the identification
of $T_h\mathcal T$ with $H^1_{Ad\rho}(\pi_1(S),\mathbb R^3)$)
is
\[
  e^\lambda_r(h)(\alpha)=\sum m_i\xi(w_i)=\xi\left(\sum m_iw_i\right)\,.
\]
Comparing this formula with  Equations (\ref{mess:eq}) and (\ref{mess:eq2}) 
we deduce that
$e^\lambda_r(h)=\xi_*(\tau_-^{\lambda}(h))$.
\end{proof}